\input ./style/arxiv-general.cfg
\documentclass[aos,MSNbibl,nameyear,dvips]{arximspdf}
\makeatletter
   \@ifpackageloaded{graphicx}{}{\usepackage{graphicx}}
\makeatother
\usepackage{mathbh}
\usepackage[ruled]{algorithm2e}

\doi{10.1214/15-AOS1369}
\volume{44}
\issue{5}
\pubyear{2016}
\firstpage{1896}
\lastpage{1930}
\docsubty{FLA}

\makeatletter
\let\my@algocf@latexcaption\algocf@latexcaption
\let\my@addcontentsline\addcontentsline
\long\def\algocf@latexcaption#1[#2]#3{%
\def\addcontentsline##1##2##3{}%
\my@algocf@latexcaption{#1}[#2]{#3}%
\global\let\addcontentsline\my@addcontentsline%
}
\def\sfrac#1#2{#1/#2}
\def\vfrac#1#2{(#1)/#2}

\def\vafrac#1#2{(#1)/(#2)}
\def\sklfrac#1#2{(#1/#2)}

\newcommand{\rrvert}{\vert}
\newcommand{\rrVert}{\Vert}
\newcommand{\llvert}{\vert}
\newcommand{\llVert}{\Vert}
\newcommand{\mmid}{\|}
\renewcommand{\mid}{|}

\newcommand{\eqref}[1]{(\ref{#1})}
\newcommand{\mathbbm}{\mathbh}
\newtheorem{teo}{Theorem}
\newtheorem{prop}[teo]{Proposition}
\newtheorem{lemma}[teo]{Lemma}
\newcommand{\argmin}{\mathop{\operatorname{argmin}}}
\newcommand{\argmax}{\mathop{\operatorname{argmax}}}
\newcommand{\sargmin}{\mathop{\operatorname{sargmin}}}
\newcommand{\sargmax}{\mathop{\operatorname{sargmax}}}
\makeatother

\begin{document}
\begin{frontmatter}

\title{Statistical and computational trade-offs in estimation of
sparse principal components}
\runtitle{Computational bounds in sparse PCA}

\begin{aug}
\author[A]{\fnms{Tengyao}~\snm{Wang}\thanksref{M1,T1}\ead[label=e1]{t.wang@statslab.cam.ac.uk}\ead[label=u1,url]{http://www.statslab.cam.ac.uk/\textasciitilde rjs57}},
\author[A]{\fnms{Quentin}~\snm{Berthet}\thanksref{M1,M2,T2}\ead[label=e2]{q.berthet@statslab.cam.ac.uk}\ead[label=u2,url]{http://www.statslab.cam.ac.uk/\textasciitilde tw389}}\\
\and
\author[A]{\fnms{Richard~J.}~\snm{Samworth}\corref{}\thanksref{M1,T3}\ead[label=e3]{r.samworth@statslab.cam.ac.uk}\ead
[label=u3,url]{http://www.statslab.cam.ac.uk/\texttildelow qb204}}
\runauthor{T. Wang, Q. Berthet and R. J. Samworth}
\affiliation{University of Cambridge\thanksmark{M1} and California
Institute of Technology\thanksmark{M2}}
\address[A]{Statistical Laboratory\\
Wilberforce Road\\
Cambridge, CB3 0WB\\
United Kingdom\\
\printead{e1}\\
\phantom{E-mail: }\printead*{e2}\\
\phantom{E-mail: }\printead*{e3}\\
\printead{u1}\\
\phantom{URL: }\printead*{u2}\\
\phantom{URL: }\printead*{u3}}
\end{aug}
\thankstext{T1}{Supported by a Benefactors' scholarship from St. John's College, Cambridge.}
\thankstext{T2}{Supported by Air Force Office of Scientific Research (AFOSR) Grant FA9550-14-1-0098 at the Center for
the Mathematics of Information at the California Institute of Technology.}
\thankstext{T3}{Supported by Engineering and Physical Sciences Research Council Early Career Fellowship EP/J017213/1 and Leverhulme Trust Grant PLP-2014-353.}

\footnotetext{\textbf{Tribute}: Peter was a remarkable person: not only a prolific and
highly influential researcher, but also someone with a wonderful warmth and generosity
of spirit.  He was a great inspiration to so many statisticians around the world.
 We are deeply saddened that he is no longer with us, and dedicate this paper to his
  memory. Further personal reflections on Peter Hall's life and work from the third
  author can be found in Samworth (\citeyear{Sam06}).}

%
\received{\smonth{5} \syear{2015}}
%
\revised{\smonth{7} \syear{2015}}

%
\begin{abstract}
In recent years, sparse principal component analysis has emerged as an
extremely popular dimension reduction technique
for high-dimensional data. The theoretical challenge, in the simplest
case, is to estimate the leading eigenvector
of a population covariance matrix under the assumption that this
eigenvector is sparse. An impressive range of
estimators have been proposed; some of these are fast to compute, while
others are known to achieve the minimax
optimal rate over certain Gaussian or sub-Gaussian classes. In this
paper, we show that, under a
widely-believed assumption from computational complexity theory, there
is a fundamental trade-off
between statistical and computational performance in this problem. More
precisely, working with new,
larger classes satisfying a restricted covariance concentration
condition, we show that there is an
effective sample size regime in which no randomised polynomial time
algorithm can achieve the minimax
optimal rate. We also study the theoretical performance of a
(polynomial time) variant of the well-known
semidefinite relaxation estimator, revealing a subtle interplay between
statistical and computational efficiency.
\end{abstract}

%
\end{frontmatter}

\setcounter{footnote}{3}

\section{Introduction}\label{sec1}

Principal Component Analysis (PCA), which involves
projecting a sample of multivariate data onto the space spanned by the
leading eigenvectors of the sample covariance matrix, is one of the
oldest and most widely-used dimension reduction devices in statistics.
It has proved to be particularly effective when the dimension of the
data is relatively small by comparison with the sample size. However,
the work of
\citet{JohnstoneLu2009} and \citet{Paul2007} shows that PCA
breaks down in the high-dimensional settings that are frequently
encountered in many diverse modern application areas. For instance,
consider the spiked covariance model where $X_1,\ldots,X_n$ are
independent $N_p(0,\Sigma)$ random vectors, with $\Sigma= I_p +
\theta v_1 v_1^\top$ for some $\theta> 0$ and an arbitrary unit
vector $v_1 \in\mathbb{R}^p$. In this case, $v_1$ is the leading
eigenvector (principal component) of $\Sigma$, and the classical PCA
estimate would be $\hat{v}_1$, a unit-length leading eigenvector of
the sample covariance matrix $\hat{\Sigma}:= n^{-1}\sum_{i=1}^n X_i
X_i^\top$. In the high-dimensional setting where $p = p_n$ is such
that $p/n \rightarrow c \in(0,1)$, \citet{Paul2007} showed that
\[
\bigl\llvert\hat{v}_1^\top v_1\bigr\rrvert
\stackrel{\mathrm{a.s.}} {\rightarrow} \cases{
\displaystyle\sqrt{\frac{1-c/\theta^2}{1 +
c/\theta}}, &\quad if $\theta>
\sqrt{c}$,
\vspace*{3pt}\cr
0, &\quad if $\theta\leq\sqrt{c}$.}
\]
In other words, $\hat{v}_1$ is inconsistent as an estimator of $v_1$
in this asymptotic regime. This phenomenon is related to the so-called
``BBP'' transition in random matrix theory
\citep{BBP2005}.

Sparse principal component analysis was designed to remedy this
inconsistency and to give additional interpretability to the projected
data. In the simplest case, it is assumed that the leading eigenvector
$v_1$ of the population covariance matrix $\Sigma$ belongs to the
$k$-sparse unit Euclidean sphere in $\mathbb{R}^p$, given by
%
\begin{equation}
\label{EqB0k} B_0(k):= \Biggl\{u = (u_1,
\ldots,u_p)^\top\in\mathbb{R}^p: \sum
_{j=1}^p \mathbbm{1}_{\{u_j \neq0\}} \leq k, \llVert
u\rrVert_2 = 1 \Biggr\}.
\end{equation}
A remarkable number of recent papers have proposed estimators of $v_1$
in this setting, including \citet{JTU2003}, \citet
{ZHT2006}, \citet{dAEGJL2007}, \citet{JohnstoneLu2009},
\citet{WTH2009}, \citet{JNRS2010}, \citet{BJNP2013},
\citet{CMW2013}, \citet{Ma2013}, \citet{SSM2013} and
\citet{VuLei2013}.

Sparse PCA methods have gained high popularity in many diverse applied
fields where high-dimensional datasets are routinely handled. These
include computer vision for online visual tracking \citep{WLY2013} and
pattern recognition \citep{NYS2011}, signal processing for image
compression \citep{Majumdar2009} and electrocardiography feature
extraction \citep{JohnstoneLu2009}, and biomedical research for gene
expression analysis \citep{ZHT2006,CS09,PTB09,ChanHall2010}, RNA-seq
classification \citep{TPW2014} and metabolomics studies \citep
{GeneveraMaleticSavatic2011}. In these applications, sparse PCA is
employed to identify a small number of interpretable directions that
represent the data succinctly, typically as the first stage of a more
involved procedure such as classification, clustering or regression.

The success of the ultimate inferential methods in the types of
application described above depends critically on how well the
particular sparse PCA technique involved identifies the relevant
meaningful directions in the underlying population. It therefore
becomes important to understand the ways in which our ability to
estimate these directions from data depends on the characteristics of
the problem, including the sample size, dimensionality, sparsity level
and signal-to-noise ratio. Such results form a key component of any
theoretical analysis of an inference problem in which sparse PCA is
employed as a first step.

In terms of the theoretical properties of existing methods for sparse
PCA, \citet{Ma2013} was able to show that his estimator attains
the minimax rate of convergence over a certain Gaussian class of
distributions, provided that $k$ is treated as a fixed constant. Both
\citet{CMW2013} and \citet{VuLei2013} also study minimax
properties, but treat $k$ as a parameter of the problem that may vary
with the sample size $n$. In particular, for a certain class $\mathcal
{P}_p(n,k)$ of sub-Gaussian distributions and in a particular
asymptotic regime, \citet{VuLei2013} show\footnote{Here and
below, $a_n \asymp b_n$ means $0 < \liminf_{n \rightarrow\infty}
\llvert a_n/b_n\rrvert \leq\limsup_{n \rightarrow\infty} \llvert
a_n/b_n\rrvert < \infty$.} that
\[
\inf_{\hat{v}} \sup_{P \in\mathcal{P}_p(n,k)} \mathbb{E}_P
\bigl\{1 - \bigl(v_1^\top\hat{v}\bigr)^2\bigr
\} \asymp\frac{k \log p}{n},
\]
where the infimum is taken over all estimators $\hat{v}$; see also
\citet{BJNP2013}. Moreover, they show that the minimax rate is
attained by a leading $k$-sparse eigenvector of $\hat{\Sigma}$, given by
%
\begin{equation}
\label{Eqvmax} \hat{v}_{\max}^k \in\argmax_{u \in B_0(k)}
u^\top\hat{\Sigma} u.
\end{equation}

The papers cited above would appear to settle the question of sparse
principal component estimation (at least in a sub-Gaussian setting)
from the perspective of statistical theory. However, there remains an
unsettling feature, namely that neither the estimator of \citet
{CMW2013}, nor that of \citet{VuLei2013}, is computable in
polynomial time.\footnote{Since formal definitions of such notions
from computational complexity theory may be unfamiliar to many
statisticians, and to keep the paper as self-contained as possible, we
provide a brief introduction to this topic in Section~2 of the online supplementary material [\citet{WBSSupp2014}].} For instance, computing the estimator~\eqref{Eqvmax}
is an \textsf{NP}-hard problem, and the naive algorithm that searches
through all
${p\choose k}$ of the $k \times k$ principal submatrices of $\hat
{\Sigma}$ quickly becomes infeasible for even moderately large $p$ and $k$.

Given that sparse PCA methods are typically applied to massive
high-dimen\-sional datasets, it is crucial to understand the rates that
can be achieved using only computationally efficient procedures.
Specifically, in this paper, we address the question of whether it is
possible to find an estimator of $v_1$ that is computable in
(randomised) polynomial time, and that attains the minimax optimal rate
of convergence when the sparsity of $v_1$ is allowed to vary with the
sample size. Some progress in a related direction was made by
\citeauthor{BerthetRigollet2013a} (\citeyear{BerthetRigollet2013a,BerthetRigollet2013b}), who considered
the problem of testing the null hypothesis $H_0: \Sigma= I_p$ against
the alternative $H_1: v^\top\Sigma v \geq1+\theta$ for some $v \in
B_0(k)$ and $\theta> 0$. Of interest here is the minimal level $\theta
= \theta_{n,p,k}$ that ensures small asymptotic testing error. Under a
hypothesis on the computational intractability of a certain well-known
problem from theoretical computer science (the ``Planted Clique''
detection problem), Berthet and Rigollet showed that for certain
classes of distributions, there is a gap between the minimal $\theta
$-level permitting successful detection with a randomised polynomial
time test, and the corresponding $\theta$-level when arbitrary tests
are allowed.

The particular classes of distributions considered in \citeauthor{BerthetRigollet2013a}
(\citeyear{BerthetRigollet2013a,BerthetRigollet2013b}) were highly tailored to the
testing problem, and do not provide sufficient structure to study
principal component estimation. The thesis of this paper, however, is
that from the point of view of both theory and applications, it is the
estimation of sparse principal components, rather than testing for the
existence of a distinguished direction, that is the more natural and
fundamental (as well as more challenging) problem. Indeed, we observe
subtle phase transition phenomena that are absent from the hypothesis
testing problem; see Section~\ref{SecHighSS} for further details. It
is worth noting that different results for statistical and
computational trade-offs for estimation and testing were also observed
in the context of $k$-SAT formulas in \citet{FPV2013} and
\citet{Berthet2014}, respectively.

Our first contribution, in Section~\ref{454654211}, is to introduce a
new Restricted Covariance Concentration (RCC) condition that underpins
the classes of distributions $\mathcal{P}_p(n,k,\theta)$ over which
we perform the statistical and computational analyses [see~\eqref
{EqClasses} for a precise definition]. The RCC condition is satisfied
by sub-Gaussian distributions, and moreover has the advantage of being
more robust to certain mixture contaminations that turn out to be of
key importance in the statistical analysis under the computational
constraint. We show that subject to mild restrictions on the parameter values,
\[
\inf_{\hat{v}} \sup_{P \in\mathcal{P}_p(n,k,\theta)} \mathbb
{E}_P L(\hat{v},v_1) \asymp\sqrt{\frac{k \log p}{n\theta^2}},
\]
where $L(u,v):= \{1 - (u^\top v)^2\}^{1/2}$, and where no restrictions
are placed on the class of estimators $\hat{v}$. By contrast, in
Section~\ref{33335545445}, we show that a variant $\hat{v}^{\mathrm
{SDP}}$ of the semidefinite relaxation estimator of \citet
{dAEGJL2007} and \citet{BAd2010}, which is computable in
polynomial time, satisfies
\[
\sup_{P \in\mathcal{P}_p(n,k,\theta)} \mathbb{E}_P L\bigl(\hat
{v}^{\mathrm{SDP}},v_1\bigr) \leq(16\sqrt{2}+2)\sqrt{
\frac{k^2 \log
p}{n\theta^2}}.
\]
Our main result, in Section~\ref{4444554541}, is that, under a much
weaker planted clique hypothesis than that in \citeauthor{BerthetRigollet2013a} (\citeyear{BerthetRigollet2013a,BerthetRigollet2013b}), for any $\alpha\in
(0,1)$, there exists a moderate effective sample size asymptotic regime
in which every sequence $(\hat{v}^{(n)})$ of randomised polynomial
time estimators satisfies
\[
\sqrt{\frac{n\theta^2}{k^{1+\alpha} \log p}} \sup_{P \in\mathcal
{P}_p(n,k,\theta)} \mathbb{E}_P
L\bigl(\hat{v}^{(n)},v_1\bigr) \rightarrow\infty.
\]
%
This result shows that there is a fundamental trade-off between
statistical and computational efficiency in the estimation of sparse
principal components, and that there is in general no consistent
sequence of randomised polynomial time estimators in this regime.
Interestingly, in a high effective sample size regime, where even
randomised polynomial time estimators can be consistent, we are able to
show in Theorem~\ref{ThmRateOptimal} that under additional
distributional assumptions, a modified (but still polynomial time)
version of $\hat{v}^{\mathrm{SDP}}$ attains the minimax optimal rate.
Thus, the trade-off disappears for a sufficiently high effective sample
size, at least over a subset of the parameter space.

Statistical and computational trade-offs have also recently been
studied in the context of convex relaxation algorithms \citep
{ChandrasekaranJordan2013}, submatrix signal detection \citep
{MaWu2013,ChenXu2014}, sparse linear regression \citep{ZWJ2014},
community detection \citep{HWX2014} and sparse canonical correlation
analysis \citep{GMZ2014}. Given the importance of computationally
feasible algorithms with good statistical performance in today's era of
big data, it seems clear that understanding the extent of this
phenomenon in different settings will represent a key challenge for
theoreticians in the coming years.

Proofs of our main results are given in the \hyperref[appen]{Appendix}, while several
ancillary results are deferred to the online supplementary material
[\citet{WBSSupp2014}]. We end this section by introducing some
notation used throughout the paper. For a vector $u = (u_1,\ldots
,u_M)^\top\in\mathbb{R}^M$, a matrix $A = (A_{ij}) \in\mathbb
{R}^{M \times N}$ and for $q \in[1,\infty)$, we write $\llVert
u\rrVert _q:=
(\sum_{i=1}^M \llvert u_i\rrvert ^q )^{1/q}$ and $\llVert A\rrVert
_q:= (\sum_{i=1}^M \sum_{j=1}^N \llvert A_{ij}\rrvert ^q )^{1/q}$ for
their (entrywise)
$\ell_q$-norms. We also write $\llVert u\rrVert _0:= \sum_{i=1}^M
\mathbbm
{1}_{\{u_i \neq0\}}$, $\operatorname{supp}(u):= \{i: u_i\neq0\}$,
$\llVert A\rrVert
_0:=\break \sum_{i=1}^M \sum_{j=1}^N \mathbbm{1}_{\{A_{ij} \neq0\}}$ and
$\operatorname{supp}(A):= \{(i,j): A_{ij}\neq0\}$. For $S \subseteq\{
1,\ldots,M\}$ and $T \subseteq\{1,\ldots,N\}$, we write $u_S:=
(u_i:i \in S)^\top$ and write $M_{S,T}$ for the $\llvert S\rrvert
\times\llvert T\rrvert $
submatrix of $M$ obtained by extracting the rows and columns with
indices in $S$ and $T$, respectively. For positive sequences $(a_n)$ and
$(b_n)$, we write $a_n \ll b_n$ to mean $a_n/b_n \rightarrow0$.

\section{Restricted covariance concentration and minimax rate of estimation}
\label{454654211}

Let $p \geq2$ and let $\mathcal{P}$ denote the class of probability
distributions $P$ on $\mathbb{R}^p$ with $\int_{\mathbb{R}^p} x
\,dP(x) = 0$ and such that the entries of $\Sigma(P):= \int_{\mathbb
{R}^p} x x^\top \,dP(x)$ are finite. For $P \in\mathcal{P}$, write
$\lambda_1(P),\ldots,\lambda_p(P)$ for the eigenvalues of $\Sigma
(P)$, arranged in decreasing order. When $\lambda_1(P) - \lambda_2(P)
> 0$, the first principal component $v_1(P)$, that is, a unit-length
eigenvector of $\Sigma$ corresponding to the eigenvalue $\lambda
_1(P)$, is well defined up to sign. In some places below, and where it
is clear from the context, we suppress the dependence of these
quantities on $P$, or write the eigenvalues and eigenvectors as
$\lambda_1(\Sigma),\ldots,\lambda_p(\Sigma)$ and $v_1(\Sigma
),\ldots,v_p(\Sigma)$, respectively. 
Let $X_1,\ldots,X_n$ be independent and identically distributed random
vectors with distribution $P$, and form the $n \times p$ matrix
$\mathbf{X}:= (X_1,\ldots,X_n)^\top$. An \emph{estimator} of $v_1$
is a measurable function from $\mathbb{R}^{n \times p}$ to $\mathbb
{R}^p$, and we write $\mathcal{V}_{n,p}$ for the class of all such estimators.

Given unit vectors $u,v \in\mathbb{R}^p$, let $\Theta(u,v):= \cos
^{-1}(\llvert u^\top v\rrvert )$ denote the acute angle between $u$
and $v$, and
define the loss function
\[
L(u,v):= \sin\Theta(u,v) = \bigl\{1 - \bigl(u^\top v
\bigr)^2\bigr\}^{1/2} = \frac
{1}{\sqrt{2}}\bigl\llVert
uu^\top- vv^\top\bigr\rrVert_2.
\]
Note that $L(\cdot,\cdot)$ is invariant to sign changes of either of
its arguments. The \emph{directional variance} of $P$ along a unit
vector $u \in\mathbb{R}^p$ is defined to be $V(u):= \mathbb{E}\{
(u^\top X_1)^2\} = u^\top\Sigma u$. Its empirical counterpart is $\hat
V(u):= n^{-1} \sum_{i=1}^n (u^\top X_i)^2 = u^\top\hat\Sigma u$,
where $\hat{\Sigma}:= n^{-1}\sum_{i=1}^n X_i X_i^\top$ denotes the
sample covariance matrix.

Recall the definition of the $k$-sparse unit ball $B_0(k)$ from~\eqref
{EqB0k}. Given $\ell\in\{1,\ldots,p\}$ and $C \in(0,\infty)$,
we say $P$ satisfies a \emph{Restricted Covariance Concentration}
(RCC) condition with parameters $p, n, \ell$ and $C$, and write $P \in
\mathrm{RCC}_p(n,\ell,C)$,~if
%
\begin{equation}
\label{EqRCC} \mathbb{P} \biggl\{\sup_{u\in B_0(\ell)} \bigl\llvert\hat
V(u) - V(u)\bigr\rrvert\geq C\max\biggl(\sqrt{\frac{\ell\log(p/\delta)}{n}},
\frac{\ell\log
(p/\delta)}{n} \biggr) \biggr\} \leq\delta
\end{equation}
for all $\delta> 0$. It is also convenient to define
\[
\mathrm{RCC}_p(\ell,C):= \bigcap_{n=1}^\infty
\mathrm{RCC}_p(n,\ell,C) \quad\mbox{and} \quad\mathrm{RCC}_p(C)
:= \bigcap_{\ell=1}^p \mathrm{RCC}_p(
\ell,C).
\]
The RCC conditions amount to uniform Bernstein-type concentration
properties of the directional variance around its expectation along all
sparse directions. This condition turns out to be particularly
convenient in the study of convergence rates in sparse PCA, and
moreover, as we show in Proposition~\ref{PropSubgaussianRCC} below,
sub-Gaussian distributions satisfy an RCC condition for all sample
sizes $n$ and all sparsity levels~$\ell$. Recall that a mean-zero
distribution $Q$ on $\mathbb{R}^p$ is \emph{sub-Gaussian} with
parameter\footnote{Note that some authors say that distributions
satisfying this condition are sub-Gaussian with parameter $\sigma$,
rather than $\sigma^2$.} $\sigma^2 \in(0,\infty)$, written
\[
Q \in\mathrm{sub\mbox{-}Gaussian}_p\bigl(\sigma^2\bigr),
\]
if whenever $Y \sim Q$, we have $\mathbb{E}(e^{u^\top Y}) \leq
e^{\sigma^2 \llVert u\rrVert ^2/2}$ for all $u \in\mathbb{R}^p$.

\begin{prop}\label{PropSubgaussianRCC}
\textup{(i)}~For every $\sigma> 0$, we have
\[
\mathrm{sub\mbox{-}Gaussian}_p\bigl(\sigma^2\bigr) \subseteq
\mathrm{RCC}_p \biggl(16\sigma^2 \biggl(1+
\frac{9}{\log p} \biggr) \biggr).
\]

\textup{(ii)} In the special case where $P = N_p(0,\Sigma)$, we have $P \in
\mathrm{RCC}_p (8\lambda_1(P)(1+\frac{9}{\log p}) )$.
\end{prop}
%

Our convergence rate results for sparse principal component estimation
will be proved over the following classes of distributions. For $\theta
> 0$, let
%
\begin{eqnarray}
\label{EqClasses} \mathcal{P}_p(n,k,\theta)&:=& \bigl\{P \in
\mathrm{RCC}_p(n,2,1) \cap\mathrm{RCC}_p(n,2k,1):
\nonumber\\[-8pt]\\[-8pt]\nonumber
&&{}v_1(P) \in B_0(k), \lambda_1(P) -
\lambda_2(P) \geq\theta\bigr\}.
\end{eqnarray}
Observe that RCC classes have the scaling property that if the
distribution of a random vector $Y$ belongs to $\mathrm{RCC}_p(n,\ell
,C)$ and if $r > 0$, then the distribution of $rY$ belongs to $\mathrm
{RCC}_p(n,\ell,r^2 C)$. It is therefore convenient to fix $C=1$ in
both RCC classes in~\eqref{EqClasses}, so that $\theta$ becomes a
measure of the signal-to-noise level. 

For a symmetric $A \in\mathbb{R}^{p \times p}$, define $\hat
{v}_{\mathrm{max}}^k(A):= \sargmax_{u \in B_0(k)} u^\top A u$ to be
the $k$-sparse maximum eigenvector of $A$, where $\sargmax$ denotes
the smallest element of the argmax in the lexicographic ordering. [This
choice ensures that $\hat{v}_{\mathrm{max}}^k(A)$ is a measurable
function of $A$.] Theorem~\ref{ThmMinimaxUpperBound} below gives a
finite-sample minimax upper bound for estimating $v_1(P)$ over
$\mathcal{P}_{p}(n,k,\theta)$. For similar bounds over Gaussian or
sub-Gaussian classes, see \citet{CMW2013} and \citet
{VuLei2013}, who consider the more general problem of principal
subspace estimation. As well as working with a larger class of
distributions, our different proof techniques facilitate an explicit constant.
%
\begin{teo}
\label{ThmMinimaxUpperBound}
For $2k \log p \leq n$, the $k$-sparse empirical maximum eigenvector,
$\hat{v}_{\mathrm{max}}^k(\hat\Sigma)$, satisfies
\[
\sup_{P\in\mathcal{P}_{p}(n,k,\theta)} \mathbb{E}_P L \bigl(\hat
{v}_{\mathrm{max}}^k(\hat\Sigma), v_1(P) \bigr) \leq2
\sqrt{2} \biggl(1+\frac{1}{\log p} \biggr)\sqrt\frac{k\log p}{n\theta
^2} \leq7\sqrt
\frac{k\log p}{n\theta^2}.
\]
\end{teo}
A matching minimax lower bound of the same order in all parameters $k,
p, n$ and $\theta$ is given below. The proof techniques are adapted
from \citet{VuLei2013}.
%
\begin{teo}
\label{ThmMinimaxLowerBound}
Suppose that $7 \leq k \leq p^{1/2}$ and $0 < \theta\leq\frac
{1}{16(1+\sfrac{9}{\log p})}$. Then
\[
\inf_{\hat v \in\mathcal{V}_{n,p}} \sup_{P \in\mathcal
{P}_p(n,k,\theta)} \mathbb{E}_P
L \bigl(\hat{v}, v_1(P) \bigr) \geq\min\biggl\{\frac{1}{1660}
\sqrt\frac{k\log p}{n\theta^2}, \frac
{5}{18\sqrt{3}} \biggr\}.
\]
\end{teo}
We remark that the conditions in the statement of Theorem~\ref
{ThmMinimaxLowerBound} can be strengthened or weakened, with a
corresponding weakening or strengthening of the constants in the bound.
For instance, a bound of the same order in $k, p, n$ and $\theta$
could be obtained assuming only that $k \leq p^{1-\delta}$ for some
$\delta> 0$. The upper bound on $\theta$ is also not particularly
restrictive. For example, if $P = N_p(0,\sigma^2I_p + \theta e_1
e_1^\top)$, where $e_1$ is the first standard basis vector in $\mathbb
{R}^p$, then it can be shown that the condition $P \in\mathcal
{P}_p(n,k,\theta)$ requires that $\theta\leq1 - \sigma^2$.


\section{Computationally efficient estimation}
\label{33335545445}

As was mentioned in the\break \hyperref[sec1]{Introduction}, the trouble with the estimator
$\hat{v}_{\mathrm{\max}}^k(\hat{\Sigma})$ of Section~\ref
{454654211}, as well as the estimator of \citet{CMW2013}, is that
there are no known polynomial time algorithms for their computation.
In this section, we therefore study the (polynomial time) semidefinite
relaxation estimator $\hat{v}^{\mathrm{SDP}}$ defined by
Algorithm~\ref{AlgocodeSDP} below. This estimator is a variant of one
proposed by \citet{dAEGJL2007}, whose support recovery properties
were studied for a particular class of Gaussian distributions and a
known sparsity level by \citet{AminiWainwright2009}.

\begin{algorithm}[t]
\SetAlgoLined
\IncMargin{1em}
\DontPrintSemicolon
\KwIn{
$\mathbf{X} = (X_1,\ldots,X_n)^\top\in\mathbb{R}^{n \times p}$,
$\lambda> 0, \varepsilon> 0$
}
\Begin{
\textbf{Step 1:} Set $\hat{\Sigma} \leftarrow n^{-1}\mathbf{X}^\top
\mathbf{X}$.\\
\textbf{Step 2:} For $f(M):= \operatorname{tr}(\hat{\Sigma}M) - \lambda
\llVert M\rrVert _1$, let $\hat M^{\varepsilon}$ be an $\varepsilon
$-maximiser of $f$
in $\mathcal{M}_1$. In other words, $\hat{M}^{\varepsilon}$ satisfies
$f(\hat{M}^{\varepsilon}) \geq\max_{M\in\mathcal{M}_1} f(M) -
\varepsilon$.\\
\textbf{Step 3:} Let $\hat{v}^{\mathrm{SDP}}:= \hat{v}^{\mathrm
{SDP}}_{\lambda,\varepsilon} \in\argmax_{u:\llVert u\rrVert _2=1}
u^\top\hat
{M}^{\varepsilon} u$.
}
\KwOut{$\hat{v}^{\mathrm{SDP}}$}
\caption{Pseudo-code for computing the semidefinite relaxation
estimator $\hat{v}^{\mathrm{SDP}}$}
\label{AlgocodeSDP}
\end{algorithm}

To motivate the main step (Step 2) of Algorithm~\ref{AlgocodeSDP}, it
is convenient to let~$\mathcal{M}$ denote the class of $p \times p$
nonnegative definite real, symmetric matrices, and let $\mathcal{M}_1
:= \{M \in\mathcal{M}: \operatorname{tr}(M) = 1\}$. Let $\mathcal
{M}_{1,1}(k^2):= \{M \in\mathcal{M}_1: \operatorname{rank}(M) = 1,
\llVert M\rrVert
_0 = k^2\}$ and observe that
\[
\max_{u \in B_0(k)} u^\top\hat{\Sigma} u = \max
_{u \in B_0(k)} \operatorname{tr}\bigl(\hat{\Sigma}uu^\top\bigr) =
\max_{M \in\mathcal
{M}_{1,1}(k^2)} \operatorname{tr}(\hat{\Sigma}M).
\]
In the final expression, the rank and sparsity constraints are
nonconvex. We therefore adopt the standard semidefinite relaxation
approach of dropping the rank constraint and replacing the sparsity
$(\ell_0)$ constraint with an $\ell_1$ penalty to obtain the convex
optimisation problem
%
\begin{equation}
\max_{M\in\mathcal{M}_1} \bigl\{\operatorname{tr}(\hat\Sigma M) -
\lambda\llVert
M\rrVert_1 \bigr\}. \label{EqSDP}
\end{equation}

We now discuss the complexity of computing $\hat{v}^{\mathrm{SDP}}$
in detail. One possible way of implementing Step~2 is to use a generic
interior-point method. However, as shown in \citet{Nesterov2005},
\citet{Nemirovski2004} and \citet{BAd2010}, certain
first-order algorithms [i.e., methods requiring $O(1/\varepsilon)$ steps
to find a feasible point achieving an $\varepsilon$-approximation of the
optimal objective function value] can significantly outperform such
generic interior-point solvers. The key idea in both \citet
{Nesterov2005} and \citet{Nemirovski2004} is that the
optimisation problem in Step 2 can be rewritten in a saddlepoint formulation:
\[
\max_{M\in\mathcal{M}_1} \operatorname{tr}(\hat{\Sigma}M) - \lambda
\llVert M
\rrVert_1 = \max_{M\in\mathcal{M}_1}\min_{U\in\mathcal{U}}
\operatorname{tr} \bigl((\hat{\Sigma}+U)M \bigr),
\]
where $\mathcal{U}:= \{U \in\mathbb{R}^{p\times p}: U^\top= U, \llVert
U\rrVert _\infty\leq\lambda\}$. The fact that $\operatorname{tr} ((\hat
{\Sigma}+U)M )$ is linear in both $M$ and $U$ makes the problem
amenable to proximal methods. In Algorithm~\ref{AlgoStep2} above, we
state a possible implementation of Step 2 of Algorithm~\ref
{AlgocodeSDP}, derived from the ``basic implementation'' in \citet
{Nemirovski2004}. In the algorithm, the $\llVert \cdot\rrVert
_2$-norm projection
$\Pi_{\mathcal{U}}(A)$ of a symmetric matrix $A = (A_{ij}) \in
\mathbb{R}^{p \times p}$ onto $\mathcal{U}$ is given by
\[
\bigl(\Pi_{\mathcal{U}}(A) \bigr)_{ij}:= \operatorname{sign}(A_{ij})
\min\bigl(\llvert A_{ij}\rrvert, \lambda\bigr).
\]
For the projection $\Pi_{\mathcal{M}_1}(A)$, first decompose $A =:
\mathit{PDP}^\top$ for some orthogonal~$P$ and diagonal $D = \operatorname
{diag}(d)$, where $d = (d_1,\ldots,d_p)^\top\in\mathbb{R}^p$. Now
let $\Pi_{\mathcal{W}}(d)$ be the projection\vspace*{1pt} image of $d$ on the unit
$(p-1)$-simplex $\mathcal{W}:= \{(w_1,\ldots,w_p): w_j \geq0, \sum
_{j=1}^p w_j = 1\}$. Finally,\vspace*{1pt} transform back to obtain\break $\Pi_{\mathcal
{M}_1}(A):= P\operatorname{diag} (\Pi_{\mathcal{W}}(d) )P^\top
$. The fact that Algorithm~\ref{AlgoStep2} outputs an $\varepsilon
$-maximiser of the optimisation problem in Step 2 of Algorithm~\ref
{AlgocodeSDP} follows from \citeauthor{Nemirovski2004} [(\citeyear{Nemirovski2004}), Theorem~3.2],
which implies in our particular case that after $N$ iterations,
%
\[
\max_{M\in\mathcal{M}_1} \min_{U\in\mathcal{U}} \operatorname{tr} \bigl((
\hat{\Sigma}+U)M \bigr) - \min_{U\in\mathcal{U}} \operatorname{tr} \bigl
((\hat{
\Sigma}+U)\hat{M}^\varepsilon\bigr) \leq\frac
{\lambda^2p^2+1}{\sqrt{2}N}.
\]
%

\begin{algorithm}[t]
\SetAlgoLined
\IncMargin{1em}
\DontPrintSemicolon
\KwIn{
$\hat{\Sigma} \in\mathcal{M}$, $\lambda> 0$, $\varepsilon> 0$.
}
\Begin{
Set $M_0 \leftarrow I_p/p$, $U_0 \leftarrow0 \in\mathbb{R}^{p \times
p}$ and $N \leftarrow\lceil\frac{\lambda^2p^2+1}{\sqrt
{2}\varepsilon} \rceil$.\\
\For{$t \leftarrow1$ \KwTo$N$}{
$U'_t \leftarrow\Pi_{\mathcal{U}} (U_{t-1} - \frac{1}{\sqrt
{2}}M_{t-1} ), M_t' \leftarrow\Pi_{\mathcal{M}_1}
(M_{t-1} + \frac{1}{\sqrt{2}}\hat{\Sigma} + \frac{1}{\sqrt
{2}}U_{t-1} )$.\\\vspace*{1pt}
$U_t \leftarrow\Pi_{\mathcal{U}} (U_{t-1} - \frac{1}{\sqrt
{2}}M'_t ), M_t \leftarrow\Pi_{\mathcal{M}_1} (M_{t-1} +
\frac{1}{\sqrt{2}}\hat{\Sigma} + \frac{1}{\sqrt{2}}U'_t )$.
}
Set $\hat{M}^\varepsilon\leftarrow\frac{1}{N}\sum_{t=1}^N M_t'$.
}
\KwOut{$\hat{M}^\varepsilon$}

\caption{A possible implementation of Step 2 of Algorithm~\protect
\ref{AlgocodeSDP}}
\label{AlgoStep2}
\end{algorithm}

In Algorithm~\ref{AlgocodeSDP}, Step~1 takes $O(np^2)$ floating point
operations; Step~3 takes $O(p^3)$ operations in the worst case, though
other methods such as the Lanczos method \citep
{Lanczos1950,GolubVanLoan1996} require only $O(p^2)$ operations under
certain conditions. Our particular implementation (Algorithm~\ref
{AlgoStep2}) for Step 2 requires $O(\frac{\lambda^2p^2 + 1}{\varepsilon
})$ iterations in the worst case, though this number may often be
considerably reduced by terminating the \textbf{for} loop if the
primal-dual gap
\[
\lambda_1(\hat{U}_t + \hat{\Sigma}) - \bigl\{
\operatorname{tr}(\hat{M}_t\hat{\Sigma}) - \lambda\llVert
\hat{M}_t\rrVert_1\bigr\}
\]
falls below $\varepsilon$, where $\hat{U}_t:= t^{-1}\sum_{s=1}^t U_s'$
and $\hat{M}_t:= t^{-1}\sum_{s=1}^t M_s'$. The most costly step
within the \textbf{for} loop is the eigen-decomposition used to
compute the projection $\Pi_{\mathcal{M}_1}$, which takes $O(p^3)$
operations. Taking $\lambda:= 4\sqrt{\frac{\log p}{n}}$ and
$\varepsilon:= \frac{\log p}{4n}$ as in Theorem~\ref{ThmRateSDP}
below, we find an overall complexity for the algorithm of $O (\max
(p^5,\frac{np^3}{\log p}) )$ operations in the worst case.

We now turn to the theoretical properties of the estimator $\hat
{v}^{\mathrm{SDP}}$ computed using Algorithm~\ref{AlgocodeSDP}.
Lemma~\ref{ThmInterRateSDP} below is stated in a general,
deterministic fashion, but will be used in Theorem~\ref{ThmRateSDP}
below to bound the loss incurred by the estimator on the event that the
sample and population covariance matrices are close in $\ell_\infty
$-norm. See also \citeauthor{VCLR2013} [(\citeyear{VCLR2013}), Theorem~3.1] for a closely
related result in the context of a projection matrix estimation
problem. Recall that $\mathcal{M}$ denotes the class of $p \times p$
nonnegative definite real, symmetric matrices.
%
\begin{lemma}
\label{ThmInterRateSDP}
Let $\Sigma\in\mathcal{M}$ be such that $\theta:= \lambda_1(\Sigma
) - \lambda_2(\Sigma) > 0$. Let $\mathbf{X} \in\mathbb{R}^{n
\times p}$ and $\hat{\Sigma}:= n^{-1}\mathbf{X}^\top\mathbf{X}$.
For arbitrary $\lambda> 0$ and $\varepsilon> 0$, if $\llVert \hat
{\Sigma} -
\Sigma\rrVert _\infty\leq\lambda$, then the semidefinite relaxation
estimator $\hat{v}^{\mathrm{SDP}}$ in Algorithm~\ref{AlgocodeSDP}
with inputs $\mathbf{X}, \lambda, \varepsilon$ satisfies
\[
L \bigl(\hat{v}^{\mathrm{SDP}},v_1(\Sigma) \bigr) \leq
\frac{4\sqrt
{2}\lambda k}{\theta} + 2\sqrt\frac{\varepsilon}{\theta}.
\]
\end{lemma}
Theorem~\ref{ThmRateSDP} below describes the statistical properties
of the estimator $\hat{v}^{\mathrm{SDP}}$ over $\mathcal
{P}_p(n,k,\theta)$ classes. It reveals in particular that we incur a
loss of statistical efficiency of a factor of $\sqrt{k}$ compared with
the minimax upper bound in Theorem~\ref{ThmMinimaxUpperBound} in
Section~\ref{454654211} above. As well as applying Lemma~\ref
{ThmInterRateSDP} on the event $\{\llVert \hat{\Sigma} - \Sigma
\rrVert _\infty
\leq\lambda\}$, the proof relies on Lemma~5
in the online supplementary material [\citet{WBSSupp2014}], which
relates the event $\{\llVert \hat{\Sigma} - \Sigma\rrVert _\infty>
\lambda\}$
to the $\mathrm{RCC}_p(n,2,1)$ condition. Indeed, this explains why we
incorporated this condition into the definition of the $\mathcal
{P}_p(n,k,\theta)$ classes.
%
\begin{teo}
\label{ThmRateSDP}
For an arbitrary $P \in\mathcal{P}_p(n,k,\theta)$ and $X_1,\ldots
,X_n \stackrel{\mathrm{i.i.d.}}{\sim} P$, we write $\hat{v}^{\mathrm
{SDP}}(\mathbf{X})$ for the output of Algorithm~\ref{AlgocodeSDP}
with input $\mathbf{X}:= (X_1,\ldots,X_n)^\top$, $\lambda:= 4\sqrt
{\frac{\log p}{n}}$ and $\varepsilon:= \frac{\log p}{4n}$. If $4\log
p\leq n \leq k^2p^2\theta^{-2}\log p$ and $\theta\in(0,k]$, then
%
\begin{equation}
\sup_{P \in\mathcal{P}_p(n,k,\theta)} \mathbb{E}_P L \bigl(\hat
{v}^{\mathrm{SDP}}(\mathbf{X}),v_1(P) \bigr) \leq\min\biggl\{ (16
\sqrt{2} + 2)\sqrt{\frac{k^2 \log p}{n\theta^2}},1 \biggr\}. \label{EqRateSDP}
\end{equation}
\end{teo}
We remark that $\hat{v}^{\mathrm{SDP}}$ has the attractive property
of being fully adaptive in the sense that it can be computed without
knowledge of the sparsity level $k$. On the other hand, $\hat
{v}^{\mathrm{SDP}}$ is not necessarily $k$-sparse. If a specific
sparsity level is desired in a particular application, Algorithm~\ref
{AlgocodeSDP} can be modified to obtain a (nonadaptive) $k$-sparse
estimator having similar estimation risk. Specifically, we can find
\[
\hat{v}^{\mathrm{SDP}}_0 \in\argmin_{u \in B_0(k)} L\bigl(\hat
{v}^{\mathrm{SDP}},u\bigr).
\]
Since $L(\hat{v}^{\mathrm{SDP}}, u)^2 = 1- (u^\top\hat
{v}^{\mathrm{SDP}} )^2$, we can compute $\hat{v}^{\mathrm
{SDP}}_0$ by setting all but the top $k$ coordinates of $\hat
{v}^{\mathrm{SDP}}$ in absolute value to zero and renormalising the
vector. In particular, $\hat{v}^{\mathrm{SDP}}_0$ is computable in
polynomial time. We deduce that under the same conditions as in
Theorem~\ref{ThmRateSDP}, for any $P \in\mathcal{P}_p(n,k,\theta)$,
\begin{eqnarray*}
&&\mathbb{E} L \bigl(\hat{v}^{\mathrm{SDP}}_0, v_1
\bigr)
\\
&&\qquad \leq\mathbb{E} \bigl[ \bigl\{L \bigl(\hat{v}^{\mathrm{SDP}}_0,
\hat{v}^{\mathrm{SDP}} \bigr) + L \bigl(\hat{v}^{\mathrm{SDP}}, v_1
\bigr) \bigr\}\mathbbm{1}_{\{\llVert \hat{\Sigma}-\Sigma\rrVert
_\infty
\leq\lambda\}} \bigr] + \mathbb{P} \bigl(\llVert\hat{
\Sigma} - \Sigma\rrVert_\infty> \lambda\bigr)
\\
&&\qquad \leq2\mathbb{E} \bigl\{L \bigl(\hat{v}^{\mathrm{SDP}}_0,
v_1 \bigr)\mathbbm{1}_{\{\llVert \hat{\Sigma}-\Sigma\rrVert _\infty
\leq\lambda\}
} \bigr\} + \mathbb{P} \bigl(
\llVert\hat{\Sigma} - \Sigma\rrVert_\infty> \lambda\bigr)
\\
&&\qquad \leq(32\sqrt{2}+3)\sqrt\frac{k^2\log p}{n\theta^2},
\end{eqnarray*}
where the final inequality follows from the proof of Theorem~\ref{ThmRateSDP}.


\section{Computational lower bounds in sparse principal component estimation}
\label{4444554541}

Theorems~\ref{ThmRateSDP} and~\ref{ThmMinimaxUpperBound} reveal a
gap between the provable performance of our semidefinite relaxation
estimator $\hat{v}^{\mathrm{SDP}}$ and the minimax optimal rate. It
is natural to ask whether there exists a computationally efficient
algorithm that achieves the statistically optimal rate of convergence.
In fact, as we will see in Theorem~\ref{ThmCompLowerBound} below, the
effective sample size region over which $\hat{v}^{\mathrm{SDP}}$ is
consistent is essentially tight among the class of all \emph
{randomised polynomial time algorithms.}\footnote{In this section,
terms from computational complexity theory defined Section~2 of the online supplementary material [\citet
{WBSSupp2014}] are written in italics at their first occurrence.}
Indeed, any randomised polynomial time algorithm with a faster rate of
convergence could otherwise be adapted to solve instances of the
planted clique problem that are believed to be hard; see Section~\ref
{SecPCproblem} below for formal definitions and discussion. In this
sense, the extra factor of $\sqrt{k}$ is an intrinsic price in
statistical efficiency that we have to pay for computational
efficiency, and the estimator $\hat{v}^{\mathrm{SDP}}$ studied in
Section~\ref{33335545445} has essentially the best possible rate of
convergence among computable estimators.
\subsection{The planted clique problem}
\label{SecPCproblem}

A \emph{graph} $G:= (V(G),E(G))$ is an ordered pair in which $V(G)$
is a countable set, and $E(G)$ is a subset of $ \{\{x,y\}:x,y \in
V(G),x \neq y \}$. For $x,y \in V(G)$, we say $x$ and $y$ are
\emph{adjacent}, and write $x \sim y$, if $\{x,y\} \in E(G)$. A \emph
{clique} $C$ is a subset of $V(G)$ such that $\{x,y\} \in E(G)$ for all
distinct $x, y \in C$. The problem of finding a clique of maximum size
in a given graph $G$ is known to be \textsf{NP}-\emph{complete}
\citep{Karp1972}. It is therefore natural to consider randomly
generated input graphs with a clique ``planted'' in, where the signal is
much less confounded by the noise. Such problems were first suggested
by \citet{Jerrum1992} and \citet{Kucera1995} as a
potentially easier variant of the classical clique problem.

Let $\mathbb{G}_m$ denote the collection of all graphs with $m$
vertices. Define $\mathcal{G}_m$ to be the distribution on $\mathbb
{G}_m$ associated with the standard Erd\H{o}s--R\'enyi random graph.
In other words, under $\mathcal{G}_m$, each pair of vertices is
adjacent independently with probability $1/2$. For any $\kappa\in\{
1,\ldots,m\}$, let $\mathcal{G}_{m,\kappa}$ be a distribution on
$\mathbb{G}_m$ constructed by first picking $\kappa$ distinct
vertices uniformly at random and connecting all edges (the ``planted
clique''), then joining each remaining pair of distinct vertices by an
edge independently with probability $1/2$. The planted clique problem
has input graphs randomly sampled from the distribution $\mathcal
{G}_{m,\kappa}$. Due to the random nature of the problem, the goal of
the planted clique problem is to find (possibly randomised) algorithms
that can locate a maximum clique $K_m$ with high probability.

It is well known that, for a standard Erd\H{o}s--R\'enyi graph, $\frac
{\llvert K_m\rrvert }{2\log_2 m} \stackrel{\mathrm{a.s.}}{\rightarrow
} 1$ [e.g., \citet{GrimmettMcDiarmid1975}]. In fact, if $\kappa= \kappa
_m$ is
such that
\[
\liminf_{m \rightarrow\infty} \frac{\kappa}{2\log_2 m}
> 1,
\]
it can be shown that the planted clique is asymptotically almost
surely also the unique maximum clique in the input graph. As observed
in \citet{Kucera1995}, there exists $C > 0$ such that, if $\kappa
> C \sqrt{m\log m}$, then asymptotically almost surely, vertices in
the planted clique have larger degrees than all other vertices, in
which case they can be located in $O(m^2)$ operations. \citet
{Alonetal1998} improved the above result by exhibiting a spectral
method that, given any $c > 0$, identifies planted cliques of size
$\kappa\geq c\sqrt{m}$ asymptotically almost surely.

Although several other polynomial time algorithms have subsequently
been discovered for the $\kappa\geq c\sqrt{m}$ case
[e.g., \citet{FeigeKrauthgamer2000,FeigeRon2010,AmesVavasis2011}],
there is
no known randomised polynomial time algorithm that can detect below
this threshold. \citet{Jerrum1992} hinted at the hardness of this
problem by showing that a specific Markov chain approach fails to work
when $\kappa= O(m^{1/2-\delta})$ for some $\delta> 0$. \citet
{FeigeKrauthgamer2003} showed that Lov\`acz--Schrijiver semidefinite
programming relaxation methods also fail in this regime. \citet
{Feldmanetal2013} recently presented further evidence of the hardness
of this problem by showing that a broad class of algorithms, which they
refer to as ``statistical algorithms'', cannot solve the planted clique
problem with $\kappa= O(m^{1/2-\delta})$ in randomised polynomial
time, for any $\delta>0$. It is now widely accepted in theoretical
computer science that the planted clique problem is hard, in the sense
that the following assumption holds with $\tau= 0$:
%
\begin{enumerate}[(A1)($\tau$)]
\item[(A1)($\tau$)] For any sequence $\kappa= \kappa_m$ such that
$\kappa\leq m^\beta$ for some $0<\beta<1/2-\tau$, there is no
randomised polynomial time algorithm that can correctly identify the
planted clique with probability tending to 1 as $m \rightarrow\infty$.
\end{enumerate}
%
We state the assumption in terms of a general parameter $\tau\in
[0,1/2)$, because it will turn out below that even if only (A1)($\tau
$) holds for some $\tau\in(0,1/6)$, there are still regimes of
$(n,p,k,\theta)$ in which no randomised polynomial time algorithm can
attain the minimax optimal rate.

Researchers have used the hardness of the planted clique problem as an
assumption to prove various impossibility results in other problems.
Examples include cryptographic applications \citep
{JuelsPeinado2000,Applebaumetal2010}, testing $k$-wise independence
\citep{Alonetal2007}
and approximating Nash equilibria \citep{HazanKrauthgamer2011}. Recent
works by \citeauthor{BerthetRigollet2013a} (\citeyear{BerthetRigollet2013a,BerthetRigollet2013b}) and
\citet{MaWu2013} used a stronger hypothesis on the hardness of
detecting the presence of a planted clique to establish computational
lower bounds in sparse principal component detection and sparse
submatrix detection problems, respectively. Our assumption~(A1)(0)
assumes only the computational intractability of identifying the entire
planted clique, so in particular, is implied by hypothesis $\mathrm
{A}_{\mathrm{PC}}$ of \citet{BerthetRigollet2013b} and
Hypothesis~1 of \citet{MaWu2013}.


\subsection{Computational lower bounds}

In this section, we use a reduction argument to show that, under
assumption~(A1)($\tau$), it is impossible to achieve the statistically
optimal rate of sparse principal component estimation using randomised
polynomial time algorithms. For $\rho\in\mathbb{N}$, and for $x \in
\mathbb{R}$, we let $[x]_\rho$ denote $x$ in its binary
representation, rounded to $\rho$ significant figures. Let $[\mathbb
{R}]_\rho:= \{[x]_\rho: x \in\mathbb{R}\}$. We say $(\hat
{v}^{(n)})$ is a \emph{sequence of randomised polynomial time
estimators} of $v_1 \in\mathbb{R}^{p_n}$ if $\hat{v}^{(n)}$ is a
measurable function from $\mathbb{R}^{n \times p_n}$ to $\mathbb
{R}^{p_n}$ and if, for every $\rho\in\mathbb{N}$, there exists a
randomised polynomial time algorithm $M_{\mathrm{pr}}$ such that for
any $\mathbf{x} \in([\mathbb{R}]_\rho)^{n \times p_n}$ we have
$[\hat{v}^{(n)}(\mathbf{x})]_\rho= [M_{\mathrm{pr}}(\mathbf
{x})]_\rho$. The sequence of semidefinite programming estimators
$(\hat{v}^{\mathrm{SDP}})$ defined in Section~\ref{33335545445} is
an example of a sequence of randomised polynomial time estimators of $v_1(P)$.
%
\begin{teo}
\label{ThmCompLowerBound}
Fix $\tau\in[0,1/6)$, assume \textup{(A1)($\tau$)}, and let $\alpha\in
(0,\frac{1-6\tau}{1-2\tau})$. For any $n \in\mathbb{N}$, let
$(p,k,\theta) = (p_n, k_n, \theta_n)$ be parameters indexed by $n$
such that $k=O(p^{1/2-\tau-\delta})$ for some $\delta\in(0,1/2-\tau
)$, $n = o(p \log p)$ and $\theta\leq k^2/(1000p)$. Suppose further that
\[
\frac{k^{1+\alpha}\log p}{n\theta^2} \rightarrow0
\]
as $n \rightarrow\infty$. Let $\mathbf{X}$ be an $n \times p$ matrix
with independent rows, each having distribution $P$. Then every
sequence $(\hat v^{(n)})$ of randomised polynomial time estimators of
$v_1(P)$ satisfies
\[
\sqrt{\frac{n\theta^2}{k^{1+\alpha}\log p}} \sup_{P \in\mathcal
{P}_p(n,k,\theta)} \mathbb{E}_P
L \bigl(\hat{v}^{(n)}(\mathbf{X}),v_1(P) \bigr)
\rightarrow\infty
\]
as $n \rightarrow\infty$.
\end{teo}
We note that the choices of parameters in the theorem imply that
%
\begin{equation}
\label{EqParamRegime2} \liminf_{n\to\infty} \frac{k^2\log p}{n\theta
^2}\geq\liminf
_{n\to\infty} \frac{p}{k^2} = \infty.
\end{equation}
As remarked in Section~\ref{SecPCproblem} above, the main interest in
this theorem comes from the case $\tau=0$. Here, our result reveals
not only that no randomised polynomial time algorithm can attain the
minimax optimal rate, but also that in the effective sample size regime
described by~\eqref{EqParamRegime2}, and provided the other side
conditions of Theorem~\ref{ThmCompLowerBound} hold, there is in
general no consistent sequence of randomised polynomial time
estimators. This is in contrast to Theorem~\ref
{ThmMinimaxUpperBound}, where we saw that consistent estimation with a
computationally inefficient procedure is possible in the asymptotic
regime~\eqref{EqParamRegime2}. A further consequence of Theorem~\ref
{ThmCompLowerBound} is that, since any sequence $(p,k,\theta) =
(p_n,k_n,\theta_n)$ satisfying the conditions of Theorem~\ref
{ThmCompLowerBound} also satisfies the conditions of Theorem~\ref
{ThmRateSDP} for large $n$, the conclusion of Theorem~\ref
{ThmRateSDP} cannot be improved in terms of the exponent of $k$ (at
least, not uniformly over the parameter range given there). As
mentioned in the \hyperref[sec1]{Introduction}, for a sufficiently large effective
sample size, where even randomised polynomial time estimators can be
consistent, the statistical and computational trade-off revealed by
Theorems~\ref{ThmMinimaxUpperBound} and~\ref{ThmCompLowerBound} may
disappear. See Section~\ref{SecHighSS} below for further details, and
\citet{GMZ2014} for recent extensions of these results to
different classes of distributions.

Even though assumption~(A1)(0) is widely believed, we also present
results under the weaker family of conditions~(A1)($\tau$) for $\tau
\in(0,1/6)$ to show that a statistical and computational trade-off
still remains for certain parameter regimes even in these settings. The
reason for assuming $\tau< 1/6$ is to guarantee that there is a regime
of parameters $(n,p,k,\theta)$ satisfying the conditions of the
theorem. Indeed, if $\tau\in[0,1/6)$ and $\alpha\in(0,\frac
{1-6\tau}{1-2\tau})$, we can set $p = n$, $k=n^{1/2 - \tau- \delta
}$ for some $\delta\in(0,\frac{1}{2}-\tau- \frac{1}{3-\alpha
} )$, $\theta= k^2/(1000n)$, and in that case,
\[
\frac{k^{1+\alpha}\log p}{n\theta^2} = \frac{10^6 n\log
n}{k^{3-\alpha}} \rightarrow0,
\]
as required.

\subsection{Sketch of the proof of Theorem~\texorpdfstring{\protect\ref{ThmCompLowerBound}}{6}}
The proof of Theorem~\ref{ThmCompLowerBound} relies on a randomised
polynomial time reduction from the planted clique problem to the sparse
principal component estimation problem. The reduction is adapted from
the ``bottom-left transformation'' of \citet{BerthetRigollet2013b},
and requires a rather different and delicate analysis.

In greater detail, suppose for a contradiction that we were given a
randomised polynomial time algorithm $\hat{v}$ for the sparse PCA
problem with a rate $\sup_{P\in\mathcal{P}_p(n,k,\theta)} \mathbb
{E}_P L(\hat{v}, v_1) \leq\sqrt\frac{k^{1+\alpha}\log p}{n\theta
^2}$ for some $\alpha< 1$. Set $m \approx p \log p$ and $\kappa
\approx k \log p$, so we are in the regime where (A1)($\tau$) holds.
Given any graph $G \sim\mathcal{G}_{m,\kappa}$ with planted clique
$K \subseteq V(G)$, we draw $n+p$ vertices $u_1,\ldots,u_n,w_1,\ldots
,w_p$ uniformly at random without replacement from $V(G)$. On average
there are about $\kappa/\log\kappa$ clique vertices in $\{w_1,\ldots
,w_p\}$, and our initial aim is to identify a large fraction of these
vertices. To do this, we form an $n\times p$ matrix $\mathbf{A}:=
(\mathbbm{1}_{u_i\sim w_j})_{i,j}$, which is an off-diagonal block of
the adjacency matrix of $G$. We then replace each 0 in $\mathbf{A}$
with $-1$ and flip the signs of each row independently with probability
$1/2$ to obtain a new matrix $\mathbf{X}$. Each component of the $i$th
row of $\mathbf{X}$ has a marginal Rademacher distribution, but if
$u_i$ is a clique vertex, then the components $\{j:w_j \in K\}$ are
perfectly correlated. Writing $\bolds{\gamma}':=(\mathbbm
{1}_{\{w_j \in K\}})_{j=1,\ldots,p}$, the leading eigenvector of
$\mathbb{E}\{\mathbf{X}^\top\mathbf{X}/n\mid\bolds{\gamma}'\}$
is proportional to $\bolds{\gamma}'$, which suggests that a
spectral method might be able to find $\{w_1,\ldots,w_p\} \cap K$ with
high probability. Unfortunately, the joint distribution of the rows of
$\mathbf{X}$ is difficult to deal with directly, but since $n$ and $p$
are small relative to $m$, we can approximate $\bolds{\gamma}'$
by a random vector $\bolds{\gamma}$ having independent $\operatorname
{Bern}(\kappa/m)$ components. We can then approximate $\mathbf{X}$ by
a matrix $\mathbf{Y}$, whose rows are independent conditional on
$\bolds{\gamma}$ and have the same marginal distribution
conditional on $\bolds{\gamma} = g$ as the rows of $\mathbf{X}$
conditional on $\bolds{\gamma}' = g$.

\begin{table}[b]
\tabcolsep=0pt
\caption{Rate of convergence of best estimator in different asymptotic regimes}\label{TabSummary}
\begin{tabular*}{\tablewidth}{@{\extracolsep{\fill}}@{}lccc@{}}
\hline
& $n\ll\frac{k\log p}{\theta^2}$ & $ \frac{k\log p}{\theta^2}\ll n \ll\frac{k^2\log p}{\theta^2}$ & $n \gg\frac{k^2\log p}{\theta^2}$
\\
\hline
All estimators & $\asymp1$ & $\asymp\sqrt\frac{k\log p}{n\theta^2}$ & $\asymp\sqrt\frac{k\log p}{n\theta ^2}$\\[10pt]
Polynomial time estimators & $\asymp1$ & $\asymp1$ & $\lesssim \sqrt\frac{k^2\log p}{n\theta^2}$\\
\hline
\end{tabular*}
\end{table}

It turns out that the distribution of an appropriately scaled version
of an arbitrary row of $\mathbf{Y}$, conditional on $\bolds
{\gamma} = g$, belongs to $\mathcal{P}_p(n,k,\theta)$ for $g$
belonging to a set of high probability. We could therefore apply our
hypothetical randomised polynomial time sparse PCA algorithm to the
scaled version of the matrix $\mathbf{Y}$ to find a good estimate of
$\bolds{\gamma}$, and since $\bolds{\gamma}$ is close to
$\bolds{\gamma}'$, this accomplishes our initial goal. With high
probability, the remaining vertices in the planted clique are those
having high connectivity to the identified clique vertices in $\{
w_1,\ldots,w_p\}$, which contradicts the hypothesis~(A1)($\tau$).

\subsection{Computationally efficient optimal estimation on
subparameter spaces in the high effective sample size regime}
\label{SecHighSS}

Theorems~\ref{ThmMinimaxUpperBound}, \ref{ThmMinimaxLowerBound},
\ref{ThmRateSDP} and~\ref{ThmCompLowerBound} enable us to
summarise, in Table~\ref{TabSummary} below, our knowledge of the best
possible rate of estimation in different asymptotic regimes, both for
arbitrary statistical procedures and for those that are computable in
randomised polynomial time. (For ease of exposition, we omit here the
additional, relatively mild, side constraints required for the above
theorems to hold.) The fact that Theorem~\ref{ThmCompLowerBound} is
primarily concerned with the setting in which $\frac{k^2\log
p}{n\theta^2} \rightarrow\infty$ raises the question of whether
computationally efficient procedures could attain a faster rate of
convergence in the high effective sample size regime where $n\gg\frac
{k^2\log p}{\theta^2}$.

The purpose of this section is to extend the ideas of \citet
{AminiWainwright2009} to show that, indeed, a variant of the estimator
$\hat{v}^{\mathrm{SDP}}$ introduced in Section~\ref{33335545445}
attains the minimax optimal rate of convergence in this asymptotic
regime, at least over a subclass of the distributions in $\mathcal
{P}_p(n,k,\theta)$. \citet{Ma2013} and \citet
{YuanZhang2013} show similar results for an iterative thresholding
algorithm for other subclasses of $\mathcal{P}_p(n,k,\theta)$ under
an extra upper bound condition on $\lambda_2(P)/\lambda_1(P)$; see
also \citet{WLL2014} and \citet{DeshpandeMontanari2014}.

Let $\mathcal{T}$ denote the set of nonnegative definite matrices
$\Sigma\in\mathbb{R}^{p \times p}$ of the form
\[
\Sigma= \theta v_1 v_1^\top+ \pmatrix{
I_k & 0
\vspace*{3pt}\cr
0 & \Gamma_{p-k}},
\]
where $v_1 \in\mathbb{R}^p$ is a unit vector such that $S:=
\operatorname{supp}(v_1)$ has cardinality $k$ and where $\Gamma_{p-k}
\in\mathbb
{R}^{(p-k)\times(p-k)}$ is nonnegative definite and satisfies
$\lambda_1(\Gamma_{p-k})\leq1$. [Here, and in the proof of
Theorem~\ref{ThmRateOptimal} below, the block matrix notation refers
to the $(S,S)$, $(S,S^c)$, $(S^c,S)$ and $(S^c,S^c)$ blocks.] We now
define a subclass of distributions
\[
\tilde{\mathcal{P}}_p(n,k,\theta):= \biggl\{P\in\mathcal
{P}_p(n,k,\theta): \Sigma(P) \in\mathcal{T}, \min
_{j\in
S}\llvert v_{1,j}\rrvert\geq16\sqrt{
\frac{k\log p}{n\theta^2}} \biggr\}.
\]
We remark that $\tilde{\mathcal{P}}_p(n,k,\theta)$ is nonempty only
if $\sqrt\frac{k^2\log p}{n\theta^2}\leq\frac{1}{16}$, since
\[
1 = \llVert v_{1,S}\rrVert_2 \geq k^{1/2}\min
_{j\in S}\llvert v_{1,j}\rrvert\geq16\sqrt
\frac
{k^2\log p}{n\theta^2}.
\]
This is one reason that the theorem below only holds in the high
effective sample size regime. Our variant of $\hat{v}^{\mathrm{SDP}}$
is described in Algorithm~\ref{AlgoVariant} below.
We remark\vspace*{1pt} that $\hat{v}^{\mathrm{MSDP}}$, like $\hat{v}^{\mathrm
{SDP}}$, is computable in polynomial time.
%
\begin{teo}
\label{ThmRateOptimal}
Assume that $X_1,\ldots, X_n\stackrel{\mathrm{i.i.d.}}{\sim} P$ for
some $P\in\tilde{\mathcal{P}}_p(n,k,\theta)$.
\begin{longlist}[(a)]
\item[(a)] Let $\lambda:= 4\sqrt\frac{\log p}{n}$. The function $f$
in Step~2 of Algorithm~\ref{AlgoVariant} has a maximiser $\hat{M}
\in\mathcal{M}_{1,1}(k^2)$ satisfying $\operatorname{sgn}(\hat{M}) =
\operatorname{sgn}(v_1 v_1^\top)$.

\begin{algorithm}[t]
\SetAlgoLined
\IncMargin{1em}
\DontPrintSemicolon
\KwIn{
$\mathbf{X} = (X_1,\ldots,X_n)^\top\in\mathbb{R}^{n \times p}$,
$\lambda> 0$, $\varepsilon> 0$, $\tau> 0$.
}
\Begin{
\textbf{Step 1:} Set $\hat{\Sigma} \leftarrow n^{-1}\mathbf{X}^\top \mathbf{X}$.\\
\textbf{Step 2:} For $f(M):= \operatorname{tr}(\hat{\Sigma}M) - \lambda \llVert M\rrVert _1$, let $\hat{M}^\varepsilon$ be an $\varepsilon $-maximiser of $f$ in $\mathcal{M}_1$.\\
\textbf{Step 3:} Let $\hat{S} \leftarrow\{j \in\{1,\ldots,p\}
: \hat M^\varepsilon_{jj} \geq\tau\}$ and $\hat{v}^{\mathrm
{MSDP}} \in\mathbb{R}^p$ by $\hat{v}_{\hat S^c}^{\mathrm{MSDP}}
\leftarrow0$ and $\hat{v}_{\hat S}^{\mathrm{MSDP}} \in\argmax
_{u\in\mathbb{R}^{\llvert \hat S\rrvert }} u^\top\hat\Sigma_{\hat
S\hat S} u$.
}
\KwOut{$\hat{v}^{\mathrm{MSDP}}$}
\caption{Pseudo-code for computing the modified semidefinite
relaxation estimator $\hat{v}^{\mathrm{MSDP}}$}
\label{AlgoVariant}
\end{algorithm}

\item[(b)] Assume that $\log p\leq n$, $\theta^2 \leq Bk^{1/2}$ for
some $B \geq1$ and $p\geq\theta(n/k)^{1/2}$. We write $\hat
{v}^{\mathrm{MSDP}}$ for the output of Algorithm~\ref{AlgoVariant}
with input parameters $\mathbf{X}:= (X_1,\ldots,X_n)^\top\in
\mathbb{R}^{n \times p}$, $\lambda:= 4\sqrt\frac{\log p}{n}$,
$\varepsilon:= (\frac{\log p}{Bn})^{5/2}$ and $\tau:= (\frac{\log
p}{Bn})^2$. Then
\[
\sup_{P\in\tilde{\mathcal{P}}_p(n,k,\theta)} \mathbb{E}_P \bigl\{ L\bigl(
\hat{v}^{\mathrm{MSDP}}, v_1\bigr) \bigr\} \leq6\sqrt
\frac{k\log
p}{n\theta^2}.
\]
\end{longlist}
\end{teo}
Theorem~\ref{ThmRateOptimal} generalises Theorem~2 of \citet
{AminiWainwright2009} in two ways: first, we relax a Gaussianity
assumption to an RCC condition; second, the leading eigenvector\vspace*{2pt} of the
population covariance matrix is not required to have nonzero entries
equal to $\pm k^{-1/2}$.

\begin{figure}[t]

\includegraphics{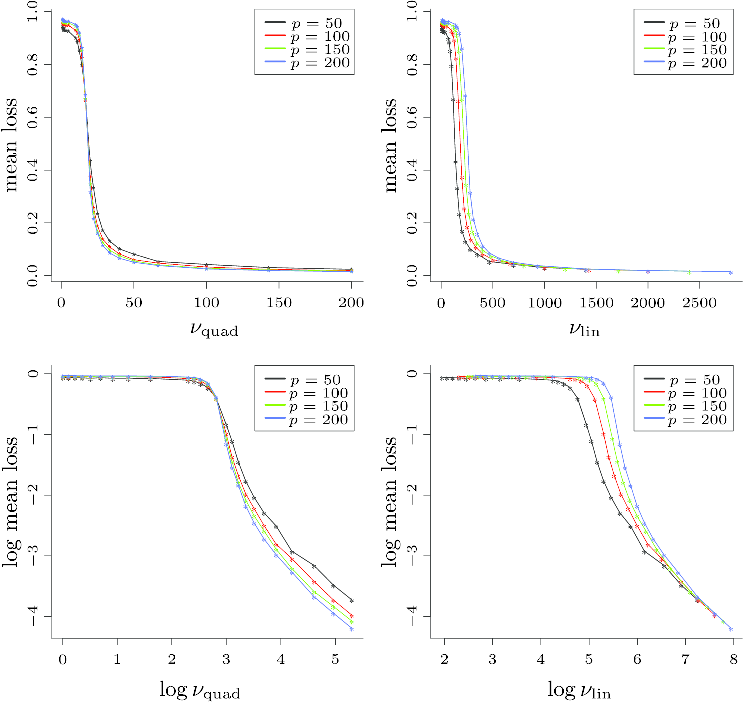}

\caption{Average loss of the estimator $\hat{v}^{\mathrm{SDP}}$ over
$N_{\mathrm{rep}} = 100$ repetitions against effective sample sizes
$\nu_{\mathrm{quad}}$ (top left) and $\nu_{\mathrm{lin}}$ (top
right). The tail behaviour under both scalings is examined under
logarithmic scales in the bottom left and bottom right panels.}
\label{FigSimulation}
\end{figure}
%
\section{Numerical experiments}
In this section, we present the results of numerical experiments to
illustrate the results of Theorems~\ref{ThmRateSDP},~\ref
{ThmCompLowerBound} and~\ref{ThmRateOptimal}. We generate $v_1 \in
\mathbb{R}^p$ by setting $v_{1,j}:= k^{-1/2}$ for $j =1,\ldots,k$,
and $v_{1,j}:= 0$ for $j = k+1,\ldots,p$. We then draw $X_1,\ldots
,X_n \stackrel{\mathrm{i.i.d.}}{\sim} N_p(0,\Sigma)$, where $\Sigma:=
I_p + \theta v_1 v_1^\top$ and $\theta=1$. We apply Algorithm~\ref
{AlgocodeSDP} to the data matrix $\mathbf{X}:= (X_1,\ldots
,X_n)^\top$ and report the average loss of the estimator $\hat
{v}^{\mathrm{SDP}}$ over $N_{\mathrm{rep}}:= 100$ repetitions. For
$p \in\{50,100,150,200\}$ and $k = \lfloor p^{1/2}\rfloor$, we repeat
the experiment for several choices of $n$ to explore the three
parameter regimes described in Table~\ref{TabSummary}. Since the
boundaries of these regimes are $n\asymp\frac{k\log p}{\theta^2}$
and $n\asymp\frac{k^2\log p}{\theta^2}$, we plot the average loss of
the experiments against effective samples sizes
\[
\nu_{\mathrm{lin}}:= \frac{n\theta^2}{k\log p}\quad\mbox{and}\quad
\nu_{\mathrm{quad}}:= \frac{n\theta^2}{k^2\log p}.
\]
The results are shown in Figure~\ref{FigSimulation}. The top left
panel of Figure~\ref{FigSimulation} shows a sharp phase transition
for the average loss, as predicted by Theorems~\ref{ThmRateSDP}
and~\ref{ThmCompLowerBound}. The right panels of Figure~\ref
{FigSimulation} suggest that in the high effective sample size regime,
$\hat{v}^{\mathrm{SDP}}$ converges at rate $\sqrt{\frac{k \log p}{n
\theta^2}}$ in this setting. This is the same rate as was proved for
the modified semidefinite relaxation estimator $\hat{v}^{\mathrm
{MSDP}}$ in Theorem~\ref{ThmRateOptimal}.

It is worth noting that it is relatively time-consuming to carry out
the simulations for the settings in the right-hand tails of the plots
in Figure~\ref{FigSimulation}. These extreme settings were chosen,
however, to illustrate that the linear scaling is the correct one in
this tail. For example, when $\nu_{\mathrm{quad}} = 200$ and $p=200$,
we require $n ={}$207,694, and the pre-processing of the data matrix to
obtain the sample covariance matrix is the time-limiting step. In
general, in our experience, the semi-definite programming algorithm is
certainly not as fast as simpler methods such as diagonal thresholding,
but is not prohibitively slow.


\begin{appendix}\label{appen}
\section{Proofs from Section~\texorpdfstring{\protect\ref{454654211}}{2}}\vspace*{-9pt}\label{AppendixGraphVector}

\begin{pf*}{Proof of Proposition~\ref{PropSubgaussianRCC}}
(i) Let $P \in\mathrm{sub\mbox{-}Gaussian}_p(\sigma^2)$, and assume that
$X_1,\ldots,X_n \stackrel{\mathrm{i.i.d.}}{\sim} P$. Then, for any $u \in
B_0(\ell)$ and $t \geq0$, we have
\[
\mathbb{P}\bigl(u^\top X_1 \geq t\bigr) \leq
e^{-t^2/\sigma^2}\mathbb{E}\bigl(e^{t
u^\top X_1/\sigma^2}\bigr) \leq e^{-t^2/(2\sigma^2)}.
\]
%
Similarly, $\mathbb{P}(-u^\top X_1 \geq t) \leq e^{-t^2/(2\sigma
^2)}$. Write $\mu_u:= \mathbb{E}\{(u^\top X_1)^2\}$; since
\begin{eqnarray*}
1 + \tfrac{1}{2}\mu_u t^2 + o
\bigl(t^2\bigr) &=& \mathbb{E}\bigl(e^{tu^\top X_1}\bigr) \leq
e^{t^2\sigma^2/2}
= 1 + \tfrac{1}{2}\sigma^2 t^2 + o
\bigl(t^2\bigr),
\end{eqnarray*}
as $t \rightarrow0$, we deduce that $\mu_u \leq\sigma^2$. Now, for
any integer $m \geq2$,
\begin{eqnarray*}
&& \mathbb{E} \bigl(\bigl\llvert\bigl(u^\top X_1
\bigr)^2 - \mu_u\bigr\rrvert^m \bigr)
\\
&&\qquad \leq \int_0^\infty
\mathbb{P} \bigl\{\bigl(u^\top X_1\bigr)^2 -
\mu_u \geq t^{1/m} \bigr\} \,dt + \mu_u^m
\\
&&\qquad \leq 2\int_0^\infty e^{-\vafrac{t^{1/m}+\mu_u}{2\sigma^2}} \,dt +
\mu_u^m
\\
&&\qquad= m!\bigl(2\sigma^2
\bigr)^m \biggl\{2e^{-\mu_u/(2\sigma^2)} + \frac
{1}{m!} \biggl(
\frac{\mu_u}{2\sigma^2} \biggr)^m \biggr\}
\\
&&\qquad \leq 2m!\bigl(2\sigma^2\bigr)^m,
\end{eqnarray*}
where the final inequality follows because the function $x \mapsto
2e^{-x} + x^m/m!$ is decreasing on $[0,1/2]$. This\vspace*{1pt} calculation allows
us to apply Bernstein's inequality [e.g., \citet{vandeGeer2000}, Lemma~5.7, taking $K
= 2\sigma^2, R = 4\sigma^2$ in her notation], to
deduce that for any $s \geq0$,
\[
\mathbb{P}\bigl(\bigl\llvert\hat V(u) - V(u)\bigr\rrvert\geq s\bigr)
\leq2\exp
\biggl(-\frac
{ns^2}{4\sigma^2 s + 32 \sigma^4} \biggr).
\]
It follows by Lemma~2 in Section~1 in the supplementary material [\citet
{WBSSupp2014}], taking $\varepsilon= 1/4$ in that result, that if $\eta>
0$ is such that $\ell\log(p/\eta) \leq n$, then for $C:= 8\sigma
^2$, we have
\begin{eqnarray*}
&& \mathbb{P} \biggl(\sup_{u \in B_0(\ell)} \bigl\llvert\hat V(u) - V(u)
\bigr\rrvert\geq2C\sqrt{\frac{\ell\log(p/\eta)}{n}} \biggr)
\\
&&\qquad \leq2\pi\ell^{1/2} \pmatrix{p\cr \ell} \biggl(\frac{128}{\sqrt
{255}}
\biggr)^{\ell-1}\exp\biggl(-\frac{C^2 \ell\log(p/\eta
)}{4C\sigma^2\sqrt{\vfrac{\ell\log(p/\eta)}{n}} + 32\sigma
^4} \biggr)
\\
&&\qquad \leq2\pi\ell^{1/2} \biggl(\frac{e}{\ell} \biggr)^\ell
\biggl(\frac
{128}{\sqrt{255}} \biggr)^{\ell-1}\eta^\ell\leq
e^9 \eta.
\end{eqnarray*}
Similarly, if $\ell\log(p/\eta) > n$, then
\begin{eqnarray*}
&& \mathbb{P} \biggl(\sup_{u \in B_0(\ell)} \bigl\llvert\hat V(u) - V(u)
\bigr\rrvert\geq2C\frac{\ell\log(p/\eta)}{n} \biggr)
\\
&&\qquad \leq2\pi\ell^{1/2} \pmatrix{p\cr \ell} \biggl(
\frac{128}{\sqrt{255}} \biggr)^{\ell-1}
\exp\biggl(-\frac{C^2 \ell^2
\log^2(p/\eta)}{4C\sigma^2\ell\log(p/\eta) + 32\sigma^4n} \biggr)
\leq e^9 \eta.
\end{eqnarray*}
Setting $\delta:= e^9 \eta$, we find (noting that we only need to
consider the case $\delta\in(0,1]$) that
\begin{eqnarray*}
&&\mathbb{P} \biggl\{\sup_{u \in B_0(\ell)} \bigl\llvert\hat V(u) - V(u)
\bigr\rrvert\geq16 \sigma^2 \biggl(1 + \frac{9}{\log p} \biggr)
\max\biggl(\sqrt{\frac{\ell\log(p/\delta)}{n}},\frac{\ell\log(p/\delta
)}{n} \biggr) \biggr\}
\\
&&\qquad \leq\mathbb{P} \biggl\{\sup_{u \in B_0(\ell)} \bigl\llvert\hat V(u) -
V(u)\bigr\rrvert\geq16 \sigma^2 \max\biggl(\sqrt{
\frac{\ell\log
(e^9p/\delta)}{n}},\frac{\ell\log(e^9p/\delta)}{n} \biggr) \biggr\}
\\
&&\qquad \leq
\delta.
\end{eqnarray*}

(ii) By Lemma~1 of \citet{LaurentMassart2000}, if $Y_1,\ldots
,Y_n$ are independent $\chi_1^2$ random variables, then for all $a > 0$,
\[
\mathbb{P} \Biggl(\frac{1}{n}\Biggl\llvert\sum
_{i=1}^n Y_i - n\Biggr\rrvert\geq a
\Biggr) \leq2e^{-\sklfrac{n}{2}(1+a-\sqrt{1+2a})} \leq2e^{-n
\min(\sfrac{a}{4},\sfrac{a^2}{16})}.
\]
Setting $\eta:= e^{-n \min(\sfrac{a}{4},\sfrac{a^2}{16})}$, we deduce that
\[
\mathbb{P} \Biggl\{\frac{1}{n}\Biggl\llvert\sum
_{i=1}^n Y_i - n\Biggr\rrvert\geq4\max
\biggl(\sqrt{\frac{\log(1/\eta)}{n}},\frac{\log
(1/\eta)}{n} \biggr) \Biggr\} \leq2\eta.
\]
Hence, using Lemma~2 again, and by a similar calculation
to part~(i),
\[
\mathbb{P} \biggl\{\sup_{u \in B_0(\ell)} \bigl\llvert\hat V(u) - V(u)
\bigr\rrvert\geq8\lambda_1(P)\max\biggl(\sqrt{\frac{\log(1/\eta)}{n}},
\frac
{\log(1/\eta)}{n} \biggr) \biggr\} \leq e^9p^\ell\eta.
\]
The result follows on setting $\delta:= e^9 p^\ell\eta$.
\end{pf*}

\begin{pf*}{Proof of Theorem~\ref{ThmMinimaxUpperBound}}
Fix an arbitrary $P \in\mathcal{P}_p(n,k,\theta)$. For notational
simplicity, we write $v:= v_1(P)$ and $\hat v:= \hat v_{\mathrm
{max}}^k(\hat\Sigma)$ in this proof. We now exploit the curvature
lemma of \citeauthor{VCLR2013} [(\citeyear{VCLR2013}), Lemma~3.1], which is closely related
to the Davis--Kahan $\sin\theta$ theorem \citep
{DavisKahan1970,YWS2015}. This lemma gives that
\[
\bigl\llVert\hat{v}\hat{v}^\top- v v^\top\bigr\rrVert
_2^2 \leq\frac{2}{\theta} \operatorname{tr} \bigl(\Sigma
\bigl(vv^\top- \hat v \hat v^\top\bigr) \bigr) \leq
\frac{2}{\theta} \operatorname{tr} \bigl((\Sigma- \hat\Sigma) \bigl
(vv^\top
- \hat v \hat v^\top\bigr) \bigr).
\]
When $\hat{v}\hat{v}^\top\neq vv^\top$, we have that $\frac
{vv^\top- \hat v \hat v^\top}{\llVert vv^\top- \hat v \hat v^\top
\rrVert _2}$
has rank 2, trace 0 and has nonzero entries in at most $2k$ rows and
$2k$ columns. It follows that its nonzero eigenvalues are $\pm1/\sqrt
{2}$, so it can be written as $(xx^\top- yy^\top)/\sqrt{2}$ for some
$x,y\in B_0(2k)$. Thus,
\begin{eqnarray*}
\mathbb{E}L(\hat{v},v) &=& \mathbb{E}\frac{1}{\sqrt{2}}\bigl\llVert\hat{v}
\hat{v}^\top- v v^\top\bigr\rrVert_2 \leq
\frac{1}{\theta} \mathbb{E} \operatorname{tr} \bigl((\Sigma- \hat\Sigma)
\bigl(xx^\top- yy^\top\bigr) \bigr)
\\
&\leq& \frac{2}{\theta}\mathbb{E}\sup_{u\in B_0(2k)} \bigl\llvert\hat
V(u) - V(u)\bigr\rrvert\leq2\sqrt{2} \biggl(1+\frac{1}{\log p} \biggr
)\sqrt
\frac
{k\log p}{n\theta^2},
\end{eqnarray*}
where we have used Proposition~1 in
Section~1 in the online supplementary material
[\citet{WBSSupp2014}] to obtain the final inequality.
\end{pf*}

\begin{pf*}{Proof of Theorem~\ref{ThmMinimaxLowerBound}}
Set $\sigma^2:= \frac{1}{8(1+\sfrac{9}{\log p})} - \theta$. We have
by Proposition~\ref{PropSubgaussianRCC}(ii) that $N_p(0, \sigma^2
I_p + \theta v_1v_1^\top) \in\mathcal{P}_p(n,k,\theta)$ for any
unit vector $v_1 \in B_0(k)$. Define $k_0:= k-1$ and $p_0:= p-1$.
Applying the variant of the Gilbert--Varshamov lemma given as
Lemma~3 in Section~1 in the online supplementary material [\citet
{WBSSupp2014}] with $\alpha:= 1/2$ and $\beta:= 1/4$, we can construct
a set $\mathcal{N}_0$ of $k_0$-sparse vectors in $\{0,1\}^{p_0}$ with
cardinality at least $(p_0/k_0)^{k_0/8}$, such that the Hamming
distance between every pair of distinct points in $\mathcal{N}_0$ is
at least $k_0$. For $\varepsilon\in(0,1]$ to be chosen later, define a
set of $k$-sparse vectors in~$\mathbb{R}^p$~by
\[
\mathcal{N}:= \Biggl\{ \pmatrix{\sqrt{1-\varepsilon^2}
\vspace*{3pt}\cr
k_0^{-1/2}\varepsilon u_0}: u_0 \in
\mathcal{N}_0 \Biggr\}.
\]
Observe that if $u,v$ are distinct elements of $\mathcal{N}$, then
\[
L(u,v) = \bigl\{1 - \bigl(u^\top v\bigr)^2\bigr
\}^{1/2} \geq\bigl\{1 - \bigl(1-\varepsilon^2/2
\bigr)^2\bigr\} ^{1/2} \geq\frac{\sqrt{3}\varepsilon}{2},
\]
and similarly $L(u,v) \leq\varepsilon$. For $u \in\mathcal{N}$, let
$P_u$ denote the multivariate normal distribution $N_p(0, \sigma^2I_p
+ \theta u u^\top)$. For any estimator $\hat v \in\mathcal
{V}_{n,p}$, we define $\hat{\psi}_{\hat{v}}:= \sargmin_{u \in
\mathcal{N}} L(\hat{v}, u)$, where $\sargmin$ denotes the smallest
element of the $\argmin$ in the lexicographic ordering. Note that $\{
\hat{\psi}_{\hat{v}} \neq u\} \subseteq\{L(\hat{v},u) \geq\sqrt
{3}\varepsilon/4\}$. We now apply the generalised version of Fano's lemma
given as\break Lemma~4 in Section~1 in the online supplementary material [\citet
{WBSSupp2014}]. Writing $D(P \mmid Q)$ for the Kullback--Leibler divergence
between two probability measures defined on the same space (a formal
definition is given just prior to Lemma~4),
we have
%
\begin{eqnarray}
\label{EqLB1}
&& \inf_{\hat v \in\mathcal{V}_{n,p}} \sup_{P \in\mathcal
{P}_p(n,k,\theta)}
\mathbb{E}_P L \bigl(\hat{v}, v_1(P) \bigr) \nonumber
\\
&&\qquad \geq \inf
_{\hat v \in\mathcal{V}_{n,p}} \max_{u \in\mathcal{N}} \mathbb{E}_{P_u}
L(\hat{v},u)
\geq \frac{\sqrt{3}\varepsilon}{4} \inf_{\hat v \in\mathcal
{V}_{n,p}} \max_{u \in\mathcal{N}}P_u^{\otimes n}(
\hat{\psi}_{\hat
{v}} \neq u)
\\
&&\qquad \geq \frac{\sqrt{3}\varepsilon}{4} \biggl(1 - \frac{\max_{u,v \in
\mathcal{N},u \neq v} D(P_v^{\otimes n} \mmid P_u^{\otimes n}) + \log
2}{(k_0/8) \log(p_0/k_0)} \biggr).\nonumber
\end{eqnarray}
We can compute, for distinct points $u,v \in\mathcal{N}$,
%
\begin{eqnarray}
\label{EqLB2} D\bigl(P_v^{\otimes n} \mmid P_u^{\otimes n}
\bigr) &=& n D(P_v \mmid P_u) = \frac{n}{2}
\operatorname{tr} \bigl( \bigl(\sigma^2 I_p + \theta u
u^\top\bigr)^{-1}\bigl(\sigma^2 I_p +
\theta v v^\top\bigr) - I_p \bigr)
\nonumber
\\
&=& \frac{n}{2} \operatorname{tr} \bigl(\bigl(\sigma^2
I_p + \theta u u^\top\bigr)^{-1} \theta
\bigl(vv^\top- uu^\top\bigr) \bigr)
\\
&=& \frac{n\theta}{2} \operatorname{tr} \biggl( \biggl(\frac{1}{\sigma^2}
I_p - \frac{\theta}{\sigma^2(\sigma^2+\theta)}uu^\top\biggr)
\bigl(vv^\top- uu^\top\bigr) \biggr)
\nonumber
\\
&=& \frac{n\theta^2}{2\sigma^2(\sigma^2+\theta)}L^2(u,v) \leq\frac
{n\theta^2\varepsilon^2}{2\sigma^2(\sigma^2+\theta)}.\nonumber
\end{eqnarray}
Let $\varepsilon:= \min\{\sqrt{a/(3b)},1\}$, where
\[
a:= 1 - \frac{8\log2}{k_0 \log(p_0/k_0)} \quad\mbox{and} \quad b:= \frac
{4n\theta^2}{\sigma^2(\sigma^2+\theta)k_0 \log(p_0/k_0)}.
\]
Then from~\eqref{EqLB1} and~\eqref{EqLB2}, we find that
\[
\inf_{\hat v \in\mathcal{V}_{n,p}} \sup_{P \in\mathcal
{P}_p(n,k,\theta)} \mathbb{E}_P
L \bigl(\hat{v}, v_1(P) \bigr) \geq\min\biggl\{\frac{1}{1660}
\sqrt\frac{k \log p}{n\theta^2}, \frac
{5}{18\sqrt{3}} \biggr\},
\]
as required.
\end{pf*}

\section{Proofs from Section~\texorpdfstring{\protect\ref{33335545445}}{3}}\vspace*{-9pt}
\label{AppendixProofsFrom3}
\begin{pf*}{Proof of Lemma~\ref{ThmInterRateSDP}}
For convenience, we write $v:= v_1(\Sigma)$, $\hat{v}$ for $\hat
{v}^{\mathrm{SDP}}$ and $\hat M$ for $\hat M^\varepsilon$ in this proof.
We first study $vv^\top- \hat{M}$, where $\hat{M} \in\mathcal
{M}_1$ is computed in Step~2 of Algorithm~\ref{AlgocodeSDP}. By the
curvature lemma of \citeauthor{VCLR2013} [(\citeyear{VCLR2013}), Lemma~3.1],
\[
\bigl\llVert vv^\top- \hat{M}\bigr\rrVert_2^2
\leq\frac{2}{\theta} \operatorname{tr} \bigl(\Sigma\bigl(vv^\top- \hat{M}
\bigr) \bigr).
\]
Moreover, since $vv^\top\in\mathcal{M}_1$, we have the basic inequality
\[
\operatorname{tr}(\hat{\Sigma}\hat{M}) - \lambda\llVert\hat{M}\rrVert_1
\geq\operatorname{tr}\bigl(\hat{\Sigma}vv^\top\bigr) - \lambda\bigl
\llVert
vv^\top\bigr\rrVert_1 - \varepsilon.
\]
Let $S$ denote the set of indices corresponding to the nonzero
components of $v$, and recall that $\llvert S\rrvert \leq k$. Since
by hypothesis $\llVert
\hat{\Sigma} - \Sigma\rrVert _\infty\leq\lambda$, we have
\begin{eqnarray*}
\bigl\llVert vv^\top- \hat{M}\bigr\rrVert_2^2
&\leq&\frac{2}{\theta} \bigl\{\operatorname{tr} \bigl(\hat{\Sigma}
\bigl(vv^\top- \hat{M}\bigr) \bigr) + \operatorname{tr} \bigl((\Sigma-
\hat{
\Sigma}) \bigl(vv^\top- \hat{M}\bigr) \bigr) \bigr\}
\\
&\leq&\frac{2}{\theta} \bigl(\lambda\bigl\llVert vv^\top\bigr\rrVert
_1 - \lambda\llVert\hat{M}\rrVert_1 + \varepsilon+
\llVert\hat{\Sigma} - \Sigma\rrVert_\infty\bigl\llVert
vv^\top- \hat{M}\bigr\rrVert_1 \bigr)
\\
&\leq&\frac{2\lambda}{\theta} \bigl(\bigl\llVert v_Sv_S^\top
\bigr\rrVert_1 - \llVert\hat{M}_{S,S}\rrVert
_1 + \bigl\llVert v_Sv_S^\top-
\hat{M}_{S,S}\bigr\rrVert_1 \bigr) + \frac
{2\varepsilon}{\theta}
\\
&\leq&\frac{4\lambda}{\theta}\bigl\llVert v_S v_S^\top-
\hat{M}_{S,S}\bigr\rrVert_1 + \frac{2\varepsilon}{\theta}\leq
\frac{4\lambda k}{\theta}\bigl\llVert v v^\top- \hat{M}\bigr\rrVert
_2+ \frac{2\varepsilon}{\theta}.
\end{eqnarray*}
We deduce that
\[
\bigl\llVert vv^\top- \hat{M}\bigr\rrVert_2 \leq
\frac{4\lambda k}{\theta} + \sqrt\frac
{2\varepsilon}{\theta}.
\]
On the other hand,
\begin{eqnarray*}
\bigl\llVert vv^\top- \hat{M}\bigr\rrVert_2^2
&=& \operatorname{tr} \bigl(\bigl(vv^\top- \hat{M}\bigr)^2 \bigr)
= 1 - 2v^\top\hat{M} v + \operatorname{tr}\bigl(\hat{M}^2\bigr)
\\
&\geq& 1 - 2\hat{v}^\top\hat{M} \hat{v} + \operatorname{tr}\bigl(
\hat{M}^2\bigr) = \bigl\llVert\hat{v}\hat{v}^\top- \hat{M}
\bigr\rrVert_2^2.
\end{eqnarray*}
We conclude that
\begin{eqnarray*}
L(\hat{v},v) &=& \frac{1}{\sqrt{2}}\bigl\llVert\hat{v}\hat{v}^\top-
vv^\top\bigr\rrVert_2 \leq\frac{1}{\sqrt{2}}\bigl(\bigl
\llVert\hat{v}\hat{v}^\top- \hat{M}\bigr\rrVert_2 + \bigl
\llVert vv^\top- \hat{M}\bigr\rrVert_2\bigr)
\\
&\leq&\sqrt{2}\bigl\llVert vv^\top- \hat{M}\bigr\rrVert_2
\leq\frac{4\sqrt{2}\lambda
k}{\theta} + 2\sqrt\frac{\varepsilon}{\theta},
\end{eqnarray*}
as required.
\end{pf*}

\begin{pf*}{Proof of Theorem~\ref{ThmRateSDP}}
Fix $P \in\mathcal{P}_p(n,k,\theta)$. By Lemma~\ref
{ThmInterRateSDP}, and by Lemma~5 in
Section~1 of the online supplementary material
[\citet{WBSSupp2014}],
%
\begin{eqnarray} \label{EqCompUpper1}
\mathbb{E}L \bigl(\hat{v}^{\mathrm{SDP}}, v_1(P) \bigr) &=& \mathbb
{E} \bigl\{ L \bigl(\hat{v}^{\mathrm{SDP}}, v_1(P) \bigr)\mathbbm
{1}_{\{\llVert \hat\Sigma- \Sigma\rrVert _\infty\leq\lambda\}} \bigr\}\nonumber
\\
&&{}+ \mathbb{E} \bigl\{ L \bigl(\hat{v}^{\mathrm{SDP}},
v_1(P) \bigr)\mathbbm{1}_{\{\llVert \hat\Sigma- \Sigma\rrVert _\infty
> \lambda
\}} \bigr\}
\\
&\leq& \frac{4\sqrt{2}\lambda k}{\theta} + 2\sqrt\frac{\varepsilon
}{\theta} + \mathbb{P} \biggl(\sup
_{u\in B_0(2)} \bigl\llvert\hat{V}(u) - V(u)\bigr\rrvert> 2\sqrt{
\frac{\log p}{n}} \biggr).\nonumber
\end{eqnarray}
Since $P\in\mathrm{RCC}_p(n,2,1)$, we have for each $\delta> 0$ that
\[
\mathbb{P} \biggl\{\sup_{u\in B_0(2)} \bigl\llvert\hat{V}(u)-V(u)
\bigr\rrvert> \max\biggl(\sqrt\frac{2\log(p/\delta)}{n}, \frac{2\log
(p/\delta)}{n} \biggr)
\biggr\} \leq\delta.
\]
Set $\delta:= \sqrt\frac{k^2\log p}{n\theta^2}$. Since $4 \log p
\leq n$, which in particular implies $n \geq3$, we have
\[
\frac{2\log(p/\delta)}{n} \leq\frac{1}{2} + \frac{1}{n}\log\biggl(
\frac{n\theta^2}{k^2 \log p} \biggr) \leq\frac{1}{2} + \frac{\log
n}{n} -
\frac{1}{n}\log\log2 \leq1.
\]
Moreover, since $n \leq k^2p^2\theta^{-2}\log p$,
\[
2\log(p/\delta) = 2\log p + \log\biggl(\frac{n\theta^2}{k^2 \log
p} \biggr) \leq4\log p.
\]
We deduce that
%
\begin{equation}
\label{EqCompUpper2} \mathbb{P} \biggl(\sup_{u\in B_0(2)} \bigl\llvert
\hat{V}(u) - V(u)\bigr\rrvert> 2\sqrt{\frac{\log p}{n}} \biggr) \leq
\sqrt
\frac{k^2\log
p}{n\theta^2}.
\end{equation}
The desired risk bound follows from~\eqref{EqCompUpper1}, the fact
that $\theta\leq k$, and~\eqref{EqCompUpper2}.
\end{pf*}

\section{Proofs from Section~\texorpdfstring{\protect\ref{4444554541}}{4}}\vspace*{-9pt}
\label{AppCompLowerBound}

\begin{pf*}{Proof of Theorem~\ref{ThmCompLowerBound}}
Suppose, for a contradiction, that there exist an infinite subset
$\mathcal{N}$ of $\mathbb{N}$, $K_0 \in[0,\infty)$ and a sequence
$(\hat{v}^{(n)})$ of randomised polynomial time estimators of $v_1(P)$
satisfying
\[
\sup_{P \in\mathcal{P}_p(n,k,\theta)} \mathbb{E}_P L \bigl(\hat
{v}^{(n)}(\mathbf{X}),v_1(P) \bigr) \leq K_0
\sqrt{\frac{k^{1+\alpha
} \log p}{n\theta^2}}
\]
for all $n \in\mathcal{N}$. 
Let $L:= \lceil\log p_n \rceil$, let $m = m_n:= \lceil
10Lp_n/9\rceil$ and let $\kappa= \kappa_n:= Lk_n$. We claim that
Algorithm~\ref{AlgoReduction} below is a randomised polynomial time
algorithm that correctly identifies the planted clique problem on $m_n$
vertices and a planted clique of size $\kappa_n$ with probability
tending to 1 as $n \rightarrow\infty$. Since $\kappa_n =
O(m_n^{1/2-\tau-\delta}\log m_n)$, this contradicts
assumption~(A1)($\tau$). We prove the claim below.

\begin{algorithm}[t]
\SetAlgoLined
\IncMargin{1em}
\DontPrintSemicolon
\KwIn{
$m \in\mathbb{N}$, $\kappa\in\{1,\ldots,m\}$, $G \in\mathbb
{G}_m$, $L \in\mathbb{N}$
}
\Begin{
\textbf{Step 1:} Let $n \leftarrow\lfloor9m/(10L) \rfloor$,
$p\leftarrow p_n$, $k \leftarrow\lfloor\kappa/L\rfloor$. Draw
$u_1,\ldots,u_n$, $w_1,\ldots,w_p$ uniformly at random without
replacement from $V(G)$. Form $\mathbf{A} = (A_{ij}) \leftarrow
(\mathbbm{1}_{\{u_i \sim w_j\}}) \in\mathbb{R}^{n\times p}$ and $
\mathbf{X} \leftarrow\operatorname{diag}(\xi_1,\ldots,\xi_n)(2\mathbf
{A} - \mathbf{1}_{n\times p})$,
where $\xi_1,\ldots,\xi_n$ are independent Rademacher random
variables (independent of $u_1,\ldots,u_n, w_1,\ldots,w_p$), and
where every entry of $\mathbf{1}_{n \times p} \in\mathbb{R}^{n
\times p}$ is 1.\\
\textbf{Step 2:} Use the randomised estimator $\hat{v}^{(n)}$ to
compute $\hat{v}= \hat{v}^{(n)}(\mathbf{X}/\sqrt{750})$.\\
\textbf{Step 3:} Let $\hat{S} = \hat{S}(\hat{v})$ be the
lexicographically smallest $k$-subset of $\{1,\ldots,p\}$ such that
$(\hat{v}_j: j\in\hat{S})$ contains the $k$ largest coordinates of
$\hat{v}$ in absolute value.\\
\textbf{Step 4:} 
For $u \in V(G)$ and $W \subseteq V(G)$, let $\operatorname{nb}(u,W):=
\mathbbm{1}_{\{u \in W\}} + \sum_{w \in W} \mathbbm{1}_{\{u \sim w\}
}$. Set $\hat{K}:= \{u \in V(G):\operatorname{nb}(u,\{w_j:j \in\hat
{S}\}) \geq3k/4 \}$.
}
\KwOut{$\hat{K}$}
\caption{Pseudo-code for a planted clique algorithm based on a
hypothetical randomised polynomial time sparse principal component
estimation algorithm}
\label{AlgoReduction}
\end{algorithm}

Let $G \sim\mathbb{G}_{m,\kappa}$, and let $K \subseteq V(G)$ denote
the planted clique. Note that the matrix $\mathbf{A}$ defined in Step
1 of Algorithm~\ref{AlgoReduction} is the off-diagonal block of the
adjacency matrix of $G$ associated with the bipartite graph induced by
the two parts $\{u_i:i=1,\ldots,n\}$ and $\{w_j:j=1,\ldots,p\}$. Let
$\bolds{\varepsilon}' = (\varepsilon_1',\ldots,\varepsilon_n')^\top$
and $\bolds{\gamma}' = (\gamma_1',\ldots,\gamma_p')^\top$,
where $\varepsilon_i':= \mathbbm{1}_{\{u_i\in K\}}$, $\gamma_j':=
\mathbbm{1}_{\{w_j \in K\}}$, and set $S':= \{j:\gamma_j' = 1\}$.

It is convenient at this point to introduce the notion of a \emph
{graph vector distribution}. We say\vspace*{1pt} $Y$ has a $p$-variate graph vector
distribution with parameters $g = (g_1,\ldots,g_p)^\top\in\{0,1\}^p$
and $\pi_0 \in[0,1]$, and write $Y \sim\mathrm{GV}_p^g(\pi_0)$, if
we can write
\[
Y = \xi\bigl\{(1-\varepsilon)R + \varepsilon(g + \tilde{R}) \bigr\},
\]
where $\xi$, $\varepsilon$ and $R$ are independent, where $\xi$ is a
Rademacher random variable, where $\varepsilon\sim\operatorname
{Bern}(\pi
_0)$, where $R = (R_1,\ldots,R_p)^\top\in\mathbb{R}^p$ has
independent Rademacher components, and where $\tilde{R} = (\tilde
{R}_1,\ldots,\tilde{R}_p)^\top$ with $\tilde{R}_j:= (1-g_j)R_j$.

Let $(\bolds{\varepsilon},\bolds{\gamma})^\top= (\varepsilon
_1,\ldots,\varepsilon_n,\gamma_1,\ldots,\gamma_p)^\top$ be $n+p$
independent Bern($\kappa/m$) random variables. For $i=1,\ldots,n$,
let $Y_i:= \xi_i \{(1-\varepsilon_i)R_i + \varepsilon_i(\bolds
{\gamma} + \tilde{R}_i) \}$ so that, conditional on $\bolds
{\gamma}$, the random vectors $Y_1,\ldots,Y_n$ are independent, each
distributed as $\mathrm{GV}_p^{\bolds{\gamma}}(\kappa/m)$. As
shorthand, we denote this conditional distribution as $Q_{\bolds
{\gamma}}$, and write $S:= \{j:\gamma_j=1\}$. Note that by
Lemma~6 in Section~1
of the online supplementary material [\citet{WBSSupp2014}],
$Q_{\bolds{\gamma}} \in\bigcap_{\ell=1}^{\lfloor20p/(9k)
\rfloor} \mathrm{RCC}_p(\ell,750)$.

Let $\bolds{Y}:= (Y_1,\ldots,Y_n)^\top$. Recall that if $P$
and $Q$ are probability measures on a measurable space $(\mathcal
{X},\mathcal{B})$, the \emph{total variation distance} between $P$
and $Q$ is defined by
\[
d_{\mathrm{TV}}(P,Q):= \sup_{B \in\mathcal{B}}\bigl\llvert P(B) - Q(B)
\bigr\rrvert.
\]
Writing $\mathcal{L}(Z)$ for the distribution (or law) of a generic
random element $Z$, and using elementary properties of the total
variation distance given in Lemma~9 in Section~1 in the online supplementary material [\citet
{WBSSupp2014}], we have
%
\begin{eqnarray}
\label{EqTV} d_{\mathrm{TV}} \bigl(\mathcal{L}(\mathbf{X}),\mathcal{L}(
\bolds{Y}) \bigr) &\leq& d_{\mathrm{TV}} \bigl(\mathcal{L} \bigl(
\bolds{\varepsilon}',\bolds{\gamma}',(R_{ij}),(
\xi_i) \bigr),\mathcal{L} \bigl(\bolds{\varepsilon},\bolds{
\gamma},(R_{ij}),(\xi_i) \bigr) \bigr)\nonumber
\\
&=& d_{\mathrm{TV}} \bigl(\mathcal{L}\bigl(\bolds{\varepsilon
'},\bolds{\gamma}'\bigr),\mathcal{L}(
\bolds{\varepsilon},\bolds{\gamma}) \bigr)
\\
&\leq& \frac{2(n+p)}{m} \leq
\frac{9(n+p)}{5p\log p}.\nonumber
\end{eqnarray}
Here, the penultimate inequality follows from \citeauthor{DiaconisFreedman1980}
[(\citeyear{DiaconisFreedman1980}), Theorem~4]. 
In view of~(\ref{EqTV}), we initially analyse Steps 2, 3 and 4 in
Algorithm~\ref{AlgoReduction} with $\mathbf{X}$ replaced by
$\bolds{Y}$. Observe that $\mathbb{E}(Y_i\mid\bolds{\gamma})
= 0$ and, writing $\Delta:= \operatorname{diag}(\bolds{\gamma}) \in
\mathbb{R}^{p \times p}$, we have
\begin{eqnarray*}
\Sigma_{\bolds{\gamma}}&:=& \operatorname{Cov}(Y_i\mid\bolds{\gamma})
= \mathbb{E} \bigl\{(1-\varepsilon_i)R_i
R_i^\top+ \varepsilon_i(\bolds{\gamma} +
\tilde R_i) (\bolds{\gamma} + \tilde{R}_i)^\top
\mid\bolds{\gamma} \bigr\}
\\
&=& I_p + \frac{\kappa}{m}\bigl(\bolds{\gamma}\bolds{
\gamma}^\top- \Delta\bigr).
\end{eqnarray*}
Writing $N_{\bolds{\gamma}}:= \sum_{j=1}^p \gamma_j$, it
follows that the largest eigenvalue of $\Sigma_{\bolds{\gamma
}}$ is $1 + \frac{\kappa}{m}(N_{\bolds{\gamma}} - 1)$, with
corresponding eigenvector $\bolds{\gamma}/N_{\bolds{\gamma
}}^{1/2} \in B_0(N_{\bolds{\gamma}})$. The other eigenvalues are~1, with multiplicity $p-N_{\bolds{\gamma}}$, and $1-\frac
{\kappa}{m}$, with multiplicity $N_{\bolds{\gamma}} - 1$.
Hence, $\lambda_1(\Sigma_{\bolds{\gamma}}) - \lambda_2(\Sigma
_{\bolds{\gamma}}) = \frac{\kappa}{m}(N_{\bolds{\gamma
}} - 1)$. Define
\[
\Gamma_0:= \biggl\{g \in\{0,1\}^p: \biggl\llvert
N_g - \frac{p\kappa
}{m}\biggr\rrvert\leq\frac{k}{20}
\biggr\},
\]
where $N_g:= \sum_{j=1}^p g_j$. We note that by Bernstein's
inequality [e.g., \citet{ShorackWellner1986}, page~855] that
%
\begin{equation}
\label{EqBound2} \mathbb{P}(\bolds{\gamma}\in\Gamma_0)
\geq1-2e^{-k/800}.
\end{equation}
If $g\in\Gamma_0$, the conditional distribution of $Y_1/\sqrt{750}$
given $\bolds{\gamma} = g$ belongs to $\mathcal{P}_p(n,k,\theta
)$ for $\theta\leq\frac{\kappa}{750m}(N_g-1)$ and all large $n \in
\mathcal{N}$. By hypothesis, it follows that for $g\in\Gamma_0$,
\[
\mathbb{E} \bigl\{L \bigl(\hat{v}^{(n)}(\bolds{Y}/\sqrt
{750}),v_1(Q_{\bolds{\gamma}}) \bigr)\mid\bolds{\gamma}=g \bigr
\} \leq K_0\sqrt{\frac{k^{1+\alpha}\log p}{n\theta^2}}
\]
for all large $n \in\mathcal{N}$. Then by Lemma~7
in Section~1 in the online supplementary
material [\citet{WBSSupp2014}], for $\hat{S}(\cdot)$ defined in
Step~3 of Algorithm~\ref{AlgoReduction}, for $g\in\Gamma_0$, and
large $n \in\mathcal{N}$,
\begin{eqnarray*}
\mathbb{E} \bigl\{\bigl\llvert S \setminus\hat{S} \bigl(\hat{v}^{(n)}(
\bolds{Y}/\sqrt{750}) \bigr)\bigr\rrvert\mid\bolds{\gamma} = g \bigr\}
&\leq& 2 N_g \mathbb{E} \bigl\{L \bigl(\hat{v}^{(n)}(
\bolds{Y}/\sqrt{750}), v_1(Q_{\bolds{\gamma}})
\bigr)^2\mid\bolds{\gamma} = g \bigr\}
\\
&\leq& 2N_g K_0\sqrt{\frac{k^{1+\alpha}\log p}{n\theta^2}}.
\end{eqnarray*}
We deduce by Markov's inequality that for $g\in\Gamma_0$, and large
$n \in\mathcal{N}$,
%
\begin{equation}
\label{EqBound3} \mathbb{P} \bigl\{\bigl\llvert S\cap\hat{S} \bigl(\hat
{v}^{(n)}(\bolds{Y}/\sqrt{750}) \bigr)\bigr\rrvert\leq
16N_{\bolds{\gamma}}/17 \mid\bolds{\gamma} = g \bigr\} \leq34 K_0
\sqrt{\frac{k^{1+\alpha}\log p}{n\theta^2}}.
\end{equation}
Let
\begin{eqnarray*}
\Omega_{0,n} &:=& \{\bolds{\gamma} \in\Gamma_0\} \cap
\bigl\{ \bigl\llvert S\cap\hat{S} \bigl(\hat{v}^{(n)}(\bolds{Y}/
\sqrt{750}) \bigr)\bigr\rrvert> 16N_{\bolds{\gamma}}/17 \bigr\},
\\
\Omega_{0,n}' &:=& \bigl\{\bolds{
\gamma}' \in\Gamma_0\bigr\} \cap\bigl\{\bigl\llvert S
\cap\hat{S} \bigl(\hat{v}^{(n)}(\bolds{X}/\sqrt{750}) \bigr)\bigr
\rrvert> 16N_{\bolds{\gamma}'}/17 \bigr\} =: \Omega_{1,n}'
\cap\Omega_{2,n}',
\end{eqnarray*}
say, where $N_{\bolds{\gamma}'}:= \sum_{j=1}^p \gamma_j'$.
When $n \in\mathcal{N}$ is sufficiently large, we have on the event
$\Omega_{0,n}'$ that
%
\begin{equation}
\label{Eq3k4} \bigl\llvert\bigl\{j \in\hat{S} \bigl(\hat{v}^{(n)}(
\bolds{X}/\sqrt{750}) \bigr): w_j \in K \bigr\}\bigr\rrvert>
3k/4.
\end{equation}
Now set
\[
\Omega_{3,n}':= \biggl\{\operatorname{nb}\bigl(u,\bigl
\{w_j:j \in S'\bigr\}\bigr) \leq\frac
{k}{2}
\mbox{ for all } u \in V(G) \setminus K \biggr\}.
\]
Recall the definition of $\hat{K}$ from Step~4 of Algorithm~\ref
{AlgoReduction}. We claim that for sufficiently large $n \in\mathcal{N}$,
\[
\Omega_{0,n}' \cap\Omega_{3,n}'
\subseteq\{\hat{K} = K\}.
\]
To see this, note that for $n \in\mathcal{N}$ sufficiently large, on
$\Omega_{0,n}'$ we have $K \subseteq\hat{K}$ by~\eqref{Eq3k4}.
For the reverse inclusion, note that if $u \in V(G) \setminus K$, then
on $\Omega_{0,n}' \cap\Omega_{3,n}'$, we have for sufficiently large
$n \in\mathcal{N}$ that
\begin{eqnarray*}
&& \operatorname{nb} \bigl(u, \bigl\{w_j:j \in\hat{S} \bigl(\hat
{v}^{(n)}(\bolds{X}/\sqrt{750}) \bigr) \bigr\} \bigr)
\\
&&\qquad \leq\bigl\llvert\{w_j:j \in\hat{S}\} \setminus K\bigr\rrvert+
\operatorname{nb} \bigl(u,\{w_j:j \in\hat{S}\} \cap K \bigr)
\\
&&\qquad \leq\bigl\llvert\{w_j:j \in\hat{S}\} \setminus K\bigr\rrvert+
\operatorname{nb} \bigl(u,\bigl\{w_j:j \in S'\bigr\} \bigr)
< \frac{k}{4} + \frac{k}{2} = \frac{3k}{4}.
\end{eqnarray*}
This establishes our claim. We conclude that for sufficiently large $n
\in\mathcal{N}$,
%
\begin{equation}
\label{EqBound0} \mathbb{P}(\hat{K} \neq K) \leq\mathbb{P} \bigl(\bigl(
\Omega_{0,n}' \cap\Omega_{3,n}'
\bigr)^c \bigr) \leq\mathbb{P} \bigl(\bigl(\Omega
_{0,n}'\bigr)^c \bigr) + \mathbb{P} \bigl(
\Omega_{1,n}' \cap\bigl(\Omega_{3,n}'
\bigr)^c \bigr).
\end{equation}
Now by Lemma 9 in Section 1 in the online
supplementary material [\citet{WBSSupp2014}], we have
%
\begin{equation}
\label{EqBound1} \bigl\llvert\mathbb{P}\bigl(\Omega_{0,n}'
\bigr) - \mathbb{P}(\Omega_{0,n})\bigr\rrvert\leq d_{\mathrm{TV}}
\bigl(\mathcal{L}\bigl(\mathbf{X},\bolds{\gamma}'\bigr),
\mathcal{L}(\mathbf{Y},\bolds{\gamma}) \bigr) \leq\frac
{9(n+p)}{5p\log p}.
\end{equation}
Moreover, by a union bound and Hoeffding's inequality, for large $n \in
\mathcal{N}$,
%
\begin{equation}
\label{EqBound4} \mathbb{P} \bigl(\Omega_{1,n}' \cap\bigl(
\Omega_{3,n}'\bigr)^c \bigr) \leq\sum
_{g \in\Gamma_0} \mathbb{P} \bigl(\bigl(\Omega_{3,n}'
\bigr)^c \mid\bolds{\gamma} = g \bigr)\mathbb{P}(\bolds{
\gamma}=g) \leq m e^{-k/800}.
\end{equation}
We conclude by~\eqref{EqBound0},~\eqref{EqBound1},~\eqref
{EqBound2},~\eqref{EqBound3} and \eqref{EqBound4} that for large
$n \in\mathcal{N}$,
\[
\mathbb{P}(\hat{K} \neq K) \leq\frac{9(n+p)}{5p\log p} + 2e^{-k/800} + 34
K_0\sqrt{\frac{k^{1+\alpha}\log p}{n\theta^2}} + m e^{-k/800}
\rightarrow0
\]
as $n \rightarrow\infty$. This contradicts assumption~(A1)($\tau$)
and, therefore, completes the proof.
\end{pf*}

\begin{pf*}{Proof of Theorem~\ref{ThmRateOptimal}}
Setting $\delta:= p^{-1}$ in~\eqref{EqRCC}, there exist events
$\Omega_1$ and $\Omega_2$, each with probability at least $1-p^{-1}$,
such that on $\Omega_1$ and $\Omega_2$, we, respectively, have
%
\begin{eqnarray}\label{EqOmega12}
\sup_{u\in B_0(2k)}\bigl\llvert\hat V(u) - V(u)\bigr\rrvert&\leq&2
\sqrt\frac{k\log p}{n} \quad\mbox{and}\quad
\nonumber\\[-8pt]\\[-8pt]\nonumber
\sup_{u\in B_0(2)}\bigl
\llvert\hat V(u) - V(u)\bigr\rrvert&\leq&2\sqrt\frac{\log p}{n}.
\end{eqnarray}
Let $\Omega_0:= \Omega_1\cap\Omega_2$. We work on $\Omega_0$
henceforth. The main ingredient for proving both parts of the theorem
is the following weak-duality inequality:
%
\begin{eqnarray}\label{EqDuality}
\max_{M\in\mathcal{M}_1} \operatorname{tr}(\hat\Sigma M) - \lambda
\llVert M
\rrVert_1 &=& \max_{M\in\mathcal{M}_1}\min_{U\in\mathcal{U}}
\operatorname{tr} \bigl((\hat\Sigma- U)M \bigr)
\nonumber
\\
&\leq&\min_{U\in\mathcal{U}} \max_{M\in\mathcal{M}_1}\operatorname{tr}
\bigl((\hat\Sigma- U)M \bigr)
\\
&=& \min_{U\in\mathcal{U}}\lambda
_1 (\hat\Sigma- U).\nonumber
\end{eqnarray}
It is convenient to denote $\gamma:= \sqrt\frac{k^2\log p}{n\theta
^2}$, and note that
\[
\gamma\leq\frac{\sqrt{k}}{16}\min_{j\in S}\llvert
v_{1,j}\rrvert\leq\frac
{1}{16}\llVert v_{1,S}\rrVert
_2 = \frac{1}{16}.
\]

\begin{pf*}{Proof of \textup{(a)}}
From~\eqref{EqDuality}, it suffices
to exhibit a primal-dual pair $(\hat M, \hat U) \in\mathcal{M}_1
\times\mathcal{U}$, such that:
\begin{longlist}[(C1)]
\item[(C1)] $\hat M = \hat v\hat v^\top$ with $\operatorname{sgn}(\hat v)
= \operatorname{sgn}(v_1)$.
\item[(C2)] $\operatorname{tr}(\hat\Sigma\hat M) - \lambda\llVert \hat
M\rrVert _1 =
\lambda_1 (\hat\Sigma- \hat U)$.
\end{longlist}
We construct the primal-dual pair as follows. Define
\[
\hat U:= \pmatrix{ \lambda\operatorname{sgn}(v_{1,S})\operatorname
{sgn}(v_{1,S})^\top& \hat\Sigma_{SS^c} -
\Sigma_{SS^c}
\vspace*{3pt}\cr
\hat\Sigma_{S^cS} - \Sigma_{S^cS}&
\hat\Sigma_{S^cS^c} - \Sigma_{S^cS^c}}.
\]
By~\eqref{EqOmega12} and Lemma~5, we have
that $\llVert \hat\Sigma- \Sigma\rrVert _\infty\leq4\sqrt\frac{\log p}{n}
\leq\lambda$, so $U\in\mathcal{U}$. Let $w = (w_1,\ldots,w_k)$ be
a unit-length leading eigenvector of $\Sigma_{SS} - \hat U_{SS}$ such
that $w^\top v_{1,S} \geq0$. Then define $\hat v$ componentwise by
\[
\hat v_S \in\mathop{\argmax_{u\in\mathbb{R}^k, \llVert u\rrVert
_2=1}}_{{u^\top w\geq0}}
u^\top(\hat\Sigma_{SS} - \hat U_{SS} ) u, \qquad
\hat v_{S^c} = 0,
\]
and set $\hat M:= \hat v \hat v^{\top}$. Note that our choices above
ensure that $\hat M \in\mathcal{M}_1$. To verify (C1), we now show
that $\operatorname{sgn}(\hat v_S) =\operatorname{sgn}(w) =
\operatorname{sgn}(v_{1,S})$. By a variant of the Davis--Kahan theorem
[\citet{YWS2015}, Theorem~2],
%
\begin{eqnarray}\label{Eqwvhat}
\llVert w - \hat v_S\rrVert_\infty& \leq& \llVert w - \hat
v_S\rrVert_2 \leq\sqrt{2} L(\hat v_S, w)
\leq\frac{2\sqrt{2}\llVert \hat\Sigma_{SS} - \Sigma_{SS}\rrVert
_{\mathrm{op}}}{\theta}
\nonumber\\[-8pt]\\[-8pt]\nonumber
&\leq&\frac{2\sqrt{2}}{\theta}\sup_{u\in B_0(2k)}\bigl\llvert\hat V(u) -
V(u)\bigr\rrvert\leq4\sqrt{2}\gamma k^{-1/2},
\end{eqnarray}
where the final inequality uses \eqref{EqOmega12}. But $w$ is also a
leading eigenvector of
\[
\frac{1}{\theta}(\Sigma_{SS} - \hat U_{SS} -
I_k) = v_{1,S} v_{1,S}^{\top} - 4 \gamma
s s^{\top},
\]
where $s:= \frac{\operatorname{sgn}(v_{1,S})}{\llVert \operatorname
{sgn}(v_{1,S})\rrVert
}$. Write $s = \alpha v_{1,S} + \beta v_{\perp}$ for some $\alpha,
\beta\in\mathbb{R}$ with $\alpha^2+\beta^2 = 1$, and a unit vector
$v_{\perp} \in\mathbb{R}^k$ orthogonal to $v_{1,S}$. Then
\begin{eqnarray*}
v_{1,S}v_{1,S}^{\top} - 4 \gamma s s^{\top} &=&
\pmatrix{ v_{1,S} & v_{\perp}}\pmatrix{1-4\gamma
\alpha^2 & -4\gamma\alpha\beta
\vspace*{3pt}\cr
-4\gamma\alpha\beta& -4\gamma
\beta^2} \pmatrix{v_{1,S}^\top
\vspace*{3pt}\cr
v_{\perp}^\top}
\\
&=& \pmatrix{ v_{1,S} & v_{\perp}}\pmatrix{a_1 &
b_1
\vspace*{3pt}\cr
a_2 & b_2} \pmatrix{d_1
& 0
\vspace*{3pt}\cr
0 & d_2} \pmatrix{a_1 & a_2
\vspace*{3pt}\cr
b_1 & b_2} \pmatrix{v_{1,S}^\top
\vspace*{3pt}\cr
v_{\perp}^\top},
\end{eqnarray*}
where $d_1 \geq d_2$ and $\pmatrix{a_1 & a_2}^\top$, $\pmatrix{b_1 &
b_2}^\top$ are eigenvalues and corresponding unit-length eigenvectors
of the middle matrix on the right-hand side of the first line. Direct
computation yields that $d_1\geq1/2 > 0\geq d_2$ and
\[
\pmatrix{a_1
\vspace*{3pt}\cr
a_2} \propto\pmatrix{1-4\gamma
\alpha^2+4\gamma\beta^2+\sqrt{16\gamma\beta^2
+ (1-4\gamma)^2}
\vspace*{3pt}\cr
-8\gamma\alpha\beta}.
\]
Consequently, $w$ is a scalar multiple of
%
\begin{equation}
a_1v_{1,S} + a_2 v_{\perp} = \bigl
\{1+4\gamma+ \sqrt{16\gamma\beta^2+(1-4\gamma)^2} \bigr\}
v_{1,S} - 8\gamma\alpha s. \label{Eqentrycompare}
\end{equation}
Since
\begin{eqnarray*}
\bigl\{1+4\gamma+ \sqrt{16\gamma\beta^2 + (1-4\gamma)^2}
\bigr\} \min_{j\in S} \llvert v_{1,j}\rrvert&\geq&2 \min
_{j\in S} \llvert v_{1,j}\rrvert \geq 32\gamma
k^{-1/2}
\\
& > & 8\gamma\alpha\llVert s\rrVert_\infty,
\end{eqnarray*}
we have $\operatorname{sgn}(w) = \operatorname{sgn}(v_{1,S})$. Hence,
by \eqref
{Eqentrycompare},
%
\begin{eqnarray}\label{Eqminw}
\min_{j=1,\ldots,k} \llvert w_j\rrvert&\geq&
\frac{ \{1+4\gamma+ \sqrt
{16\gamma\beta^2 + (1-4\gamma)^2} \}\min_{j\in S} \llvert
v_{1,j}\rrvert -
8\gamma\alpha\llVert s\rrVert _\infty}{\llVert a_1 v_{1,S} + a_2
v_{\perp}\rrVert
_2}\nonumber
\\
&\geq&\frac{(32-8\alpha)\gamma k^{-1/2}}{1+4\gamma+ \sqrt{16\gamma
\beta^2+(1-4\gamma)^2}}
\\
&\geq&\frac{12\gamma k^{-1/2}}{1+4\gamma}\geq\frac
{48}{5}\gamma
k^{-1/2}.\nonumber
\end{eqnarray}
By~\eqref{Eqwvhat} and~\eqref{Eqminw}, we have $\min_{j}\llvert
w_j\rrvert > \llVert
w-\hat v_S\rrVert _\infty$. So $\operatorname{sgn}(\hat v_S) =
\operatorname{sgn}(w) =
\operatorname{sgn}(v_{1,S})$ as desired.

It remains to check condition (C2). Since $\operatorname{sgn}(\hat v_S) =
\operatorname{sgn}(v_{1,S})$, we have
\begin{eqnarray*}
\operatorname{tr}(\hat\Sigma\hat M) - \lambda\llVert\hat M\rrVert_1 &=&
\operatorname{tr}\bigl(\hat\Sigma_{SS} \hat v_S \hat
v_S^\top\bigr) - \operatorname{tr}\bigl(\hat U_{SS}
\hat v_S \hat v_S^\top\bigr)
\\
&=& \hat v_S^\top(\hat\Sigma_{SS} - \hat
U_{SS}) \hat v_S = \lambda_1(\hat
\Sigma_{SS} - \hat U_{SS}).
\end{eqnarray*}
Moreover,
\[
\hat\Sigma- \hat U = \pmatrix{\hat\Sigma_{SS} - \hat
U_{SS} & 0
\vspace*{3pt}\cr
0 & \Gamma_{p-k}}.
\]
As $\lambda_1(\Gamma_{p-k})\leq1$ by assumption, it suffices to show
that $\lambda_1(\hat\Sigma_{SS} - \hat U_{SS})\geq1$. By Weyl's
inequality [see, e.g., \citet{HornJohnson2012}, Theorem~4.3.1]
%
\begin{eqnarray}\label{Eqlambdamax}
\lambda_1(\hat\Sigma_{SS} - \hat U_{SS}) &
\geq&\lambda_1(\Sigma_{SS} - \hat U_{SS}) -
\llVert\hat\Sigma_{SS} - \Sigma_{SS}\rrVert
_{\mathrm{op}}
\nonumber
\\[-2pt]
& \geq&1 + \theta\lambda_1\bigl(v_{1,S}v_{1,S}^\top-
4\gamma s s^\top\bigr) - 2\sqrt\frac{k\log p}{n}
\\[-2pt]
&\geq& 1 +
\frac{3\theta}{8} > 1,\nonumber
\end{eqnarray}
as required.\noqed
\end{pf*}

\begin{pf*}{Proof of \textup{(b)}}
We claim first that $\hat{S} = S$.
Let $\phi^*:= f(\hat{M})$ be the optimal value of the semidefinite
programme~\eqref{EqSDP}. From~\eqref{Eqlambdamax}, we have $\phi
^*\geq1+3\theta/8$. The proof strategy here is to use dual matrices
$\hat U$ defined in part (a) and $\hat U'$ to be defined below to
respectively bound $\operatorname{tr}(\hat M^\varepsilon_{S^cS^c})$
from above
and bound $\hat M^{\varepsilon}_{rr}$ from below for each
$r\in S$. We then check that for the choice of $\varepsilon$ we have in
the theorem, the diagonal entries of $\hat M^\varepsilon$ are above the
threshold $\log p/(6n)$ precisely when they belong to the $(S,S)$-block
of the matrix.

From~\eqref{EqDuality}, and using the fact that $\operatorname
{tr}(AB)\leq
\operatorname{tr}(A)\lambda_1(B)$ for all symmetric matrices $A$ and $B$,
we have
\begin{eqnarray*}
\operatorname{tr}\bigl(\hat\Sigma\hat M^\varepsilon\bigr) - \lambda\bigl
\llVert
\hat M^\varepsilon\bigr\rrVert_1 & \leq&\operatorname{tr} \bigl((\hat
\Sigma-\hat U)\hat M^\varepsilon\bigr)
\\
& =& \operatorname{tr} \bigl((\hat\Sigma_{SS}-\hat U_{SS})\hat
M^\varepsilon_{SS} \bigr) + \operatorname{tr} \bigl(
\Sigma_{S^cS^c} \hat M^\varepsilon_{S^cS^c} \bigr)
\\
& \leq&\operatorname{tr}\bigl(\hat M^\varepsilon_{SS}\bigr)\phi^* +
\operatorname{tr}\bigl(\hat M^\varepsilon_{S^cS^c}\bigr)\lambda_1(
\Gamma_{p-k})
\\
& =& \phi^* - \operatorname{tr}\bigl(\hat M^\varepsilon_{S^cS^c}\bigr)
\bigl(
\phi^*-1\bigr) \leq\phi^* - 3\theta\operatorname{tr}\bigl(\hat
M^\varepsilon_{S^cS^c}
\bigr)/8.
\end{eqnarray*}
On the other hand, $\operatorname{tr}(\hat\Sigma\hat M^\varepsilon) -
\lambda\llVert \hat M^\varepsilon\rrVert _1\geq\phi^* - \varepsilon$.
It follows that
%
\begin{equation}
\operatorname{tr}\bigl(\hat M^\varepsilon_{S^cS^c}\bigr) \leq
\frac{8\varepsilon}{3\theta
} \leq\frac{1}{6} \biggl(\frac{\log p}{Bn}
\biggr)^2 < \tau. \label{EqOffSignalBlock}
\end{equation}

Next, fix an arbitrary $r\in S$ and define $S_0:= S\setminus\{r\}$.
Define $\hat U'$ by
\[
\hat U'_{ij}:= \cases{ \lambda\operatorname{sgn}(\hat
M_{ij}), &\quad if $i,j\in S_0$,
\vspace*{3pt}\cr
\hat
\Sigma_{ij} - \Sigma_{ij}, &\quad otherwise.}
\]
We note that on $\Omega_0$, we have $\hat U'\in\mathcal{U}$. Again
by~\eqref{EqDuality},
%
\begin{eqnarray}\label{Eqtbc}
\operatorname{tr}\bigl(\hat\Sigma\hat M^\varepsilon\bigr) - \lambda
\bigl\llVert
\hat M^\varepsilon\bigr\rrVert_1 &\leq&\operatorname{tr} \bigl(\bigl(
\hat\Sigma-\hat U'\bigr)\hat M^\varepsilon\bigr)
\nonumber
\\
& =& \operatorname{tr} \bigl((\hat\Sigma_{S_0S_0}-\hat U_{S_0S_0})\hat
M^\varepsilon_{S_0S_0} \bigr) + \mathop{\sum
_{(i,j)\in S\times S}}_{i = r~\mathrm{or}~j =r}\Sigma_{ij}\hat
M^\varepsilon_{ji}\nonumber
\\
&&{}  + \operatorname{tr} \bigl(\Sigma_{S^cS^c} \hat
M^\varepsilon_{S^cS^c} \bigr)
\\
&\leq&\operatorname{tr}\bigl(\hat M^\varepsilon_{S_0S_0}\bigr)
\lambda_1 (\hat\Sigma_{S_0S_0} - \hat U_{S_0S_0} ) +
\mathop{\sum_{(i,j)\in S\times S}}_{i = r~\mathrm{or}~j =r}
\Sigma_{ij}\hat M^\varepsilon_{ji}\nonumber
\\
&&{} + \operatorname{tr}\bigl(
\hat M^\varepsilon_{S^cS^c}\bigr)\lambda_1(
\Gamma_{p-k}).\nonumber
\end{eqnarray}
We bound the three terms of~\eqref{Eqtbc} separately. By Lemma~8 in Section~1 in the online
supplementary material [\citet{WBSSupp2014}],
\begin{eqnarray*}
&& \lambda_1(\hat\Sigma_{S_0S_0} - \hat U_{S_0S_0})
\\
&&\qquad \leq
\lambda_1(\hat\Sigma_{SS} - \hat U_{SS}) -
\bigl\{\lambda_1(\hat\Sigma_{SS} - \hat
U_{SS}) - \lambda_2(\hat\Sigma_{SS} - \hat
U_{SS}) \bigr\}\min_{j\in S}\hat v_j^2.
\end{eqnarray*}
From~\eqref{Eqwvhat} and~\eqref{Eqminw},
\[
\min_j \llvert\hat v_j\rrvert\geq\min
_j \llvert w_j\rrvert- \llVert w-\hat
v_S\rrVert_\infty\geq3.9\gamma k^{-1/2}.
\]
Also, by Weyl's inequality,
\begin{eqnarray*}
&& \lambda_1(\hat\Sigma_{SS} - \hat U_{SS}) -
\lambda_2(\hat\Sigma_{SS} - \hat U_{SS})
\\
&&\qquad  \geq\lambda_1(\Sigma_{SS} - \hat U_{SS}) -
\lambda_2(\Sigma_{SS} - \hat U_{SS}) - 2\llVert
\hat\Sigma_{SS} - \Sigma_{SS}\rrVert_{\mathrm{op}}
\\
&&\qquad \geq\theta\bigl\{\lambda_1\bigl(v_{1,S}
v_{1,S}^\top- 4\gamma ss^\top\bigr) -
\lambda_2\bigl(v_{1,S} v_{1,S}^\top-4
\gamma ss^\top\bigr) \bigr\} - 4\sqrt\frac{k\log p}{n}
\\
&&\qquad  \geq\theta\bigl(1/2-4\gamma k^{-1/2}\bigr) \geq\theta/4.
\end{eqnarray*}
It follows that
%
\begin{equation}
\lambda_1(\hat\Sigma_{S_0S_0} - \hat U_{S_0S_0})\leq
\phi^* - 3.8\gamma^2k^{-1}\theta. \label{EqBoundFirstTerm}
\end{equation}
For the second term in~\eqref{Eqtbc}, observe that
%
\begin{eqnarray}\label{EqBoundSecondTerm}
\mathop{\sum_{(i,j)\in S\times S}}_{i = r~\mathrm{or}~j =r}\Sigma
_{ij}\hat M^\varepsilon_{ij} &\leq&\bigl(1+\theta
v_{1,r}^2\bigr)\hat M^\varepsilon_{rr} + 2
\sum_{i\in S, i\neq r} \theta v_{1,i}v_{1,r}
\hat M^\varepsilon_{i,r}
\nonumber
\\
&\leq&\hat M^\varepsilon_{rr} + 2\theta\llvert v_{1,r}
\rrvert\cdot\llVert v_1\rrVert_1\sqrt{\hat
M^\varepsilon_{rr}}
\\
&\leq& \hat M^\varepsilon_{rr} + 2
\theta\sqrt{k}\sqrt{\hat M^\varepsilon_{rr}},\nonumber
\end{eqnarray}
where the penultimate inequality uses the fact that $\hat M^\varepsilon
_{ir} \leq\sqrt{\hat M^\varepsilon_{ii}\hat M^\varepsilon_{rr}} \leq
\sqrt{\hat M^\varepsilon_{rr}}$ for a nonnegative definite matrix $\hat
M^\varepsilon$. Substituting~\eqref{EqBoundFirstTerm} and~\eqref
{EqBoundSecondTerm} into~\eqref{Eqtbc},
\begin{eqnarray*}
&& \operatorname{tr}\bigl(\hat\Sigma\hat M^\varepsilon\bigr) - \lambda
\bigl\llVert
\hat M^\varepsilon\bigr\rrVert_1
\\
&&\qquad  \leq\operatorname{tr}\bigl(\hat M^{\varepsilon}_{S_0S_0}\bigr) \biggl(\phi^* - \frac{3.8\gamma
^2\theta}{k}
\biggr) + \hat M^\varepsilon_{rr} + 2\theta\sqrt{k \hat
M^\varepsilon_{rr}} + \operatorname{tr}\bigl(\hat M^\varepsilon
_{S^cS^c}\bigr)
\\
&&\qquad  \leq\phi^* - 3.8\gamma^2k^{-1}\theta\operatorname{tr}\bigl(
\hat M^\varepsilon_{S_0S_0}\bigr) + 2\theta\sqrt{k\hat
M^\varepsilon_{rr}}
\\
&&\qquad  \leq\phi^* - 3.8\gamma^2k^{-1}\theta\bigl\{1-
\operatorname{tr}\bigl(\hat M^\varepsilon_{S^cS^c}\bigr) \bigr\} +
2\theta
\bigl(\sqrt{k}+1.9\gamma^2\bigr) \sqrt{\hat
M^\varepsilon_{rr}}.
\end{eqnarray*}
By definition, $\operatorname{tr}(\hat\Sigma\hat M^\varepsilon) -
\lambda\llVert
\hat M^\varepsilon\rrVert _1\geq\phi^* - \varepsilon$, so together
with~\eqref
{EqOffSignalBlock}, we have
%
\begin{eqnarray} \label{EqSignalBlock}
\sqrt{\hat M_{rr}^\varepsilon}& \geq&\frac{3.8\gamma^2k^{-1}\theta(1
- \vafrac{8\varepsilon}{3\theta}) - \varepsilon}{2\theta(\sqrt
{k}+1.9\gamma^2)}\nonumber
\\
&\geq&\frac{1.9\gamma^2k^{-1}(1 - \vafrac{8\varepsilon
}{3\theta})}{(\sqrt{k}+\sfrac{1.9}{256})} - \frac{\varepsilon}{2\theta
}
\\
& \geq& 1.8\gamma^2k^{-3/2} \biggl(1-\frac{8\varepsilon}{3\theta}
\biggr) - \frac{\varepsilon}{2\theta}\nonumber
\\
&\geq&\frac{1.8k^{1/2}\log p}{n\theta^2} \biggl\{1-\frac{1}{6} \biggl(
\frac{\log p}{B n} \biggr)^2 \biggr\} - \frac{1}{32} \biggl(
\frac
{\log p}{B n} \biggr)^2 \nonumber
\\
&\geq& \frac{1.4\log p}{B n} >
\tau^{1/2}.\nonumber
\end{eqnarray}
From~\eqref{EqOffSignalBlock} and~\eqref{EqSignalBlock}, we
conclude that $\hat{S} = S$, as claimed.

To conclude, by \citeauthor{YWS2015} [(\citeyear{YWS2015}), Theorem 2], on $\Omega_0$,
\[
L\bigl(\hat{v}^{\mathrm{MSDP}}, v_1\bigr) = L\bigl(
\hat{v}^{\mathrm{MSDP}}_{S}, v_{1,S}\bigr) \leq
\frac{2\llVert \hat\Sigma_{SS}- \Sigma_{SS}\rrVert _{\mathrm
{op}}}{\lambda_1(\Sigma_{SS})-\lambda_2(\Sigma_{SS})}\leq4\sqrt\frac
{k\log p}{n\theta^2},
\]
where we used~\eqref{EqOmega12} and Lemma~5
in the online supplementary material [\citet{WBSSupp2014}] in the
final bound.

For the final part of the theorem, when $p \geq\theta\sqrt{n/k}$,
\begin{eqnarray*}
\sup_{P\in\tilde{\mathcal{P}}_p(n,k,\theta)} \mathbb{E}_P \bigl\{ L\bigl(
\hat{v}^{\mathrm{MSDP}}, v_1\bigr) \bigr\} &\leq&4\sqrt
\frac{k\log
p}{n\theta^2} + \mathbb{P}\bigl(\Omega_0^c\bigr)
\\
& \leq&4\sqrt\frac{k\log p}{n\theta^2} + 2p^{-1}\leq6\sqrt\frac
{k\log p}{n\theta^2},
\end{eqnarray*}
as desired.
\end{pf*}\noqed \end{pf*}
\end{appendix}

\section*{Acknowledgements}
We thank the anonymous reviewers for helpful and constructive comments on
an earlier draft.

\begin{supplement}[id=suppA]
\stitle{Supplementary material to ``Statistical and computational trade-offs in estimation of sparse principal components''}
\slink[doi]{10.1214/15-AOS1369SUPP} 
\sdatatype{.pdf}
\sfilename{aos1369\_supp.pdf}
\sdescription{Ancillary results and a brief introduction to computational complexity theory.}
\end{supplement}

%

\printaddresses
\end{document}